\documentclass[a4paper,12pt, reqno]{amsart}
\usepackage{amsthm}
\usepackage{amsmath,amssymb,latexsym,amsfonts,mathrsfs}
\usepackage{color}
\usepackage{graphicx}
\usepackage{bbm}
\usepackage{hyperref}
\usepackage{enumitem}
\usepackage{mathtools} 
\mathtoolsset{showonlyrefs=true}

\renewcommand{\hat}{\widehat}

\renewcommand{\bar}{\overline}

\newcommand{\bC}{\ensuremath{\mathbb{C}}}

\newcommand{\bE}{\ensuremath{\mathbb{E}}}

\newcommand{\bN}{\ensuremath{\mathbb{N}}}

\newcommand{\bP}{\ensuremath{\mathbb{P}}}

\newcommand{\bR}{\ensuremath{\mathbb{R}}}

\newcommand{\bZ}{\ensuremath{\mathbb{Z}}}

\newcommand{\cB}{\ensuremath{\mathcal{B}}}

\newcommand{\cF}{\ensuremath{\mathcal{F}}}

\newcommand{\cL}{\ensuremath{\mathcal{L}}}
\newcommand{\cM}{\ensuremath{\mathcal{M}}}

\newcommand{\Real}{\operatorname{Re}}
\newcommand{\Ima}{\operatorname{Im}}

\theoremstyle{plain}
\newtheorem{Thm}{Theorem}[section]

\newtheorem{Lem}[Thm]{Lemma}
\newtheorem{Prop}[Thm]{Proposition}
\newtheorem{Cor}[Thm]{Corollary}

\theoremstyle{definition}

\newtheorem{Rem}[Thm]{Remark}
\newtheorem{Ex}[Thm]{Example}
\newtheorem*{Nota}{Notation}

\numberwithin{equation}{section}

\begin{document}

\title{Large deviations related to Dynkin--Lamperti arcsine laws for last visit times of Markov processes} 	

\author{Takahiro Mori}
\address[Takahiro Mori]{Faculty of Arts and Sciences, Kyoto Institute of Technology, Kyoto 606-8585, Japan.}
\email{tmori@kit.ac.jp}

\author{Kei Noba}
\address[Kei Noba]{Department of Mathematics, Graduate School of Science, The University of Osaka Toyonaka, Osaka 560-0043, Japan.}
\email{knoba@math.sci.osaka-u.ac.jp}

\author{Toru Sera}
\address[Toru Sera]{Department of Systems Innovation, Graduate School of Engineering Science, The University of Osaka,
Toyonaka, Osaka 560-8531, Japan.}
\email{sera.toru.es@osaka-u.ac.jp}

\begin{abstract}
We establish large deviation estimates related to the Dynkin--Lamperti arcsine laws for subordinators whose Laplace exponents are either regularly varying or comparable to regularly varying functions. By applying these results to inverse local times, we derive large deviation estimates for the last visit times to the starting point of various Hunt processes satisfying suitable on-diagonal resolvent density estimates, such as one-dimensional generalized diffusion processes, one-dimensional L\'evy processes, and Brownian motion and its subordinate process on the Sierpi\'nski gasket.
\end{abstract}


\maketitle



\section{Introduction}
L\'evy's arcsine law \cite{Lev, Lev65} for one-dimensional Brownian motion is a classical theorem in probability theory.
Since its discovery, numerous generalizations  have been developed in various directions.
Among them, a particularly important one is the Dynkin--Lamperti theorem, which characterizes the limit distribution of the position of a subordinator immediately before  it crosses a fixed level $t$, as $t\to0+$ or $t\to\infty$, under the assumption that its Laplace exponent is regularly varying with index $\alpha\in[0,1]$. See Bingham, Goldie, and Teugels \cite[Theorem 8.6.3]{BinGolTeu}, Bertoin \cite[Theorem III.6]{Ber} or Kyprianou \cite[Theorem 5.16]{Kyp} for more details.
We note that the Dynkin--Lamperti theorem has also been extended and further developed in the study of infinite ergodic transformations \cite{Tha98, Tha05, ThZw, Zwe07cpt,  KocZwe,
MelTer12, Ter15, Ter16, Ser20, Ser25ETDS}.
Recently the third author also studied  higher order approximations in the Dynkin--Lamperti theorem for subordinators in \cite{Ser25ECP}.
The first purpose of the present paper is  to establish large deviation estimates related to the Dynkin--Lamperti theorem for subordinators whose Laplace exponents are either regularly varying or bounded above and below by regularly varying functions. The second purpose is to apply them to the last visit time of a certain class of Markov processes, including one-dimensional diffusion processes, L\'evy processes and fractional diffusion processes. We are  motivated by similar large deviation estimates for occupation times of one-dimensional diffusions studied by Kasahara and Yano \cite{KasYan}, and of infinite ergodic transformations in \cite{Ser25ETDS}. We also remark that Rouault, Yor, and Zani \cite{RYZ} and Cristadoro and Pozzoli \cite{CriPoz} also studied large deviation estimates related to arcsine laws.
Before proceeding to the details, we recall L\'evy's arcsine law for the last zero of one-dimensional Brownian motion.

\begin{Nota}
Let $f(x)$ and $g(x)\neq0$ be functions defined on a neighborhood of $a$.
If $\limsup_{x\to a} |f(x)/g(x)|<\infty$, then we write $f(x)=O(g(x))$ as $x\to a$. The notation $f(x)\sim g(x)$ as $x\to a$ means that $f(x)/g(x)\to 1$ as $x\to a$.
When $f$ and $g$ are non-negative, then 
the notation $f(x) \asymp g(x)$ as $x\to a$ (respectively, for all $x$) means that there exist some constants $0<C_1\leq C_2<\infty$ such that $C_1f(x)\leq g(x)\leq C_2 f(x)$ on some neighborhood of $a$ (respectively, for all $x$). 
\end{Nota}

\begin{Ex}[Brownian motion on $\bR$]
Let $B=(B_t)_{t\geq0}$ be a Brownian motion on $\bR$ with $B_0=0$. Let $G(t)$ denote the last time $B$ visits $0$ before or at time $t>0$:
\begin{align}
	G(t)=\max\{0\leq s\leq t\::\: B_s=0\},
	\quad
	t>0.
\end{align}
L\'evy's arcsine law for the last zero states that $G(t)/t$ is arcsine distributed, that is,
\begin{align}\label{Ex:BM-AL}
    \bP\bigg(\frac{G(t)}{t}\leq s\bigg)
    =
    \int_0^s \frac{du}{\pi \sqrt{u(1-u)}}
    =
    \frac{2}{\pi}\arcsin \sqrt{s},
    \quad 0\leq s\leq 1,
    \,t>0.	
\end{align}
We refer the reader to Chung \cite[(2.6)]{Chu} for more details.
In particular, for any fixed $t>0$,
\begin{align}\label{Ex:BM-LDE}
     \bP\bigg(\frac{G(t)}{t}\leq \varepsilon\bigg)=\bP\bigg(\frac{G(t)}{t}\geq 1-\varepsilon\bigg)\sim \frac{2}{\pi}\sqrt{\varepsilon},
     \quad	\text{as $\varepsilon\to0+$.}
\end{align}	
\end{Ex}

The Dynkin--Lamperti theorem is a generalization of \eqref{Ex:BM-AL}, and can be applied to a wide class of Markov processes. See for example Bertoin \cite[Sections III.3, IV.6 and VI.3]{Ber} and Bertoin \cite[Sections 3, 8 and 9]{Ber99}.
In this paper we extend the large deviation estimate  \eqref{Ex:BM-LDE} to a broader class of Markov processes.  
To illustrate our main results, let us consider the Brownian motion on the unbounded Sierpi\'nski gasket as a typical example.

\begin{Ex}[Brownian motion on the unbounded Sierpi\'nski gasket]\label{Ex:SG}

Let $a_0=0$, $a_1=1$ and $a_2=e^{i\pi/3}\in \bC$. For $j=0,1,2$, define a contraction $F_j:\bC\to \bC$ by $F_j(z)=(z+a_j)/2$ ($z\in\bC$). Then there uniquely exists a compact set $K\subset \bC$ such that $K=\bigcup_{j=0,1,2}F_j(K)$. We call the set $K$ the ($2$-dimensional) Sierpi\'nski gasket, and $\hat{K}=\bigcup_{n=0}^\infty 2^n (K\cup e^{i2\pi/3}K)$ the unbounded Sierpi\'nski gasket. It is well-known that the Hausdorff dimensions of $K$ and $\hat{K}$ are equal to $d_f=\log3/\log 2$. Denote by $m$ the normalized $d_f$-dimensional Hausdorff measure on $\hat{K}$ with $m(K)=1$.
We write $(X_t)_{t\geq0}$ for the Brownian motion on $\hat{K}$ in the sense of Barlow and Perkins \cite{BarPer}. In the following we will see that certain estimates similar to \eqref{Ex:BM-LDE} hold in this setting.
Fix $a\in \hat{K}$ arbitrarily. Let $G(t)$ denote the last time  $(X_t)$ visits the point $a$ prior to or at time $t>0$:
\begin{align}
	G(t)=\max\{0\leq s\leq t\::\: X_s=a\},
	\quad
	t>0.
\end{align}
Let $\varepsilon:(0,\infty)\to (0,1)$ be a non-random function satisfying $\varepsilon(t)\to0$ ($t\to0+$ or $\infty$)  and $\varepsilon(t)t\to\infty$  ($t\to\infty$). 
By using our main results (Theorems \ref{Thm:LDE-long}, \ref{Thm:LDE-short}, \ref{Thm:LV-LDE-long} and \ref{Thm:LV-LDE-short}), we can show that
\begin{align}
    \bP_a\bigg(\frac{G(t)}{t}\leq \varepsilon(t)\bigg)
    \asymp
    \varepsilon(t)^{1-(d_s/2)}
    \quad
    \text{and}
    \quad
     \bP_a\bigg(\frac{G(t)}{t}\geq 1- \varepsilon(t)\bigg)
     \asymp
     \varepsilon(t)^{d_s/2},
     \label{Ex:SG-LDE}
\end{align}
as $t\to0+$ or $\infty$,
where $\bP_{a}$ denotes the law of $(X_t)$ starting from $a$, and
$d_s=2\log 3/\log 5=1.365...\in(0,2)$ denotes the spectral dimension of $\hat{K}$.

\end{Ex}

See Section \ref{Sec:SG} for the detailed proof of \eqref{Ex:SG-LDE}. Here, we provide an idea for proving \eqref{Ex:SG-LDE}.
As shown in \cite[Theorems 1.5]{BarPer}, $(X_t)$ admits a jointly continuous transition probability density (or heat kernel) $p_t(x,y)=p_t(y,x)$ ($t>0$, $x,y\in \hat{K}$) with respect to $m$, and a so-called sub-Gaussian estimate holds. In particular, we have the on-diagonal estimate
\begin{align}
	 p_t(x,x)\asymp t^{-d_s/2},
	\quad
	\text{for all $(t,x)\in (0,\infty)\times \hat{K}$}.
	\label{Ex:SG-HK-diagonal}
\end{align}
In addition, \cite[Theorem 1.9]{BarPer} implies that the $\lambda$-resolvent density $r_\lambda(x,y)=r_\lambda(y,x)=\int_0^\infty e^{-\lambda t}p_t(x,y)\,dt$ ($\lambda>0$, $x,y\in \hat{K}$) satisfies the on-diagonal estimate
\begin{align}
     r_\lambda(x,x)
    \asymp
   \lambda^{-1+(d_s/2)},
	\quad
	\text{for all $(\lambda,x)\in (0,\infty)\times \hat{K}$.}
	\label{Ex:SG-PD-diagonal}	
\end{align}
We denote by $\tau=(\tau(s))$ the inverse local time at $a$. Then the level set of $(X_t)$ for $a$, which is closed, coincides with the closure of the range of $\tau$, that is, $\{t\geq0\::\: X_t=a\}
	=
	\bar{\{\tau(s)\::\:s\geq0\}}$,  $\bP_a$-a.s.
In particular,
\begin{align}
  	G(t)=\sup\{\tau(s)\::\:s\geq0,\, \tau_s\leq t\},
  	\quad
  	t>0,
  	\,\text{$\bP_a$-a.s.}
\end{align}
Under a suitable normalization, $(\tau,\bP_a)$ is a complete subordinator (which  will be recalled in Section \ref{Sec:subordinator}) with Laplace transform
\begin{align}
   \bE_a\exp(-\lambda \tau(s))=\exp\bigg(\frac{-s}{r_\lambda(a, a)}\bigg),
   \quad
   \lambda,s >0,	
\end{align}
Hence we can study the asymptotic behavior of $G(t)$ by using the inverse local time $\tau$ and the on-diagonal resolvent density estimate \eqref{Ex:SG-PD-diagonal}. Therefore establishing the desired estimate \eqref{Ex:SG-LDE} reduces to the study of a subordinator with Laplace exponent bounded above and below by a regularly varying function.
 We emphasize that our main results also apply to broad classes of Markov processes satisfying the on-diagonal estimate \eqref{Ex:SG-PD-diagonal} for sufficiently small or large $\lambda$ with $d_s\in(0,2)$, such as real-valued generalized diffusion processes, real-valued L\'evy processes, and symmetric Hunt processes satisfying sub-Gaussian or stable-like estimates with spectral dimension less than $2$. See Sections \ref{Sec:Diffusion}, \ref{Sec:Levy} and \ref{Sec:SG} for more details.


The organization of this paper is as follows. Section \ref{Sec:subordinator} is devoted to recalling basic facts about subordinators. In Section \ref{Sec:Main} we state our main results, that is, large deviation estimates for undershoots of subordinators, which are related to the Dynkin--Lamperti theorem. In Section \ref{Sec:Tauberian-type} we show Tauberian-type estimates, which are used to prove the main results in Section \ref{sec:Proof-Main}. We recall several facts about local times for general Hunt processes in Section \ref{Sec:Dual}. Our results apply to last visit times of specific classes of Hunt processes in Sections \ref{Sec:Diffusion}, \ref{Sec:Levy} and \ref{Sec:SG}.

\section{Bernstein Functions, Subordinators and Dynkin--Lamperti Theorem}\label{Sec:subordinator}

In this section, we review some standard facts related to subordinators, and set up notations and terminologies used in subsequent sections.

Let us recall the notions of completely monotone functions and Bernstein functions, following Schilling, Song, and Vondra{\v c}ek \cite{SchSonVon}.
We denote by $\bN=\{0,1,2,\dots\}$ the set of non-negative integers. 
A function $f:(0,\infty)\to[0,\infty)$ is called \emph{completely monotone} if $f$ is infinitely differentiable and $(-d/d\lambda)^n f(\lambda)\geq0$ for all $n\in \bN$ and $\lambda>0$.
 By Bernstein's theorem \cite[Theorem 1.4]{SchSonVon}, a function $f:(0,\infty)\to[0,\infty)$ is completely monotone if and only if it is represented in the form
 \begin{align}
     f(\lambda)
     =
     \int_{[0,\infty)}
     e^{-\lambda x}\,\mu(dx),
     \quad
    \lambda>0,
    \label{CM}	
 \end{align}
 where $\mu$ denotes a measure on $[0,\infty)$ with $\int_{[0,\infty)}
     e^{-\lambda x}\,\mu(dx)<\infty$ for any $\lambda>0$. The correspondence between $f$ and $\mu$ is one-to-one. 
 
A function $\Phi:(0,\infty)\to[0,\infty)$ is called a \emph{Bernstein function} if $\Phi$ is infinitely  differentiable and its derivative $\Phi'$ is completely monotone, in other words, $(-d/d\lambda)^n \Phi(\lambda)\leq 0$ for all $n\in\bN\setminus\{0\}$ and $\lambda>0$. 
By \cite[Theorem 3.2]{SchSonVon}, a function $\Phi:(0,\infty)\to[0,\infty)$ is a Bernstein function if and only if it is represented in the form
\begin{align}
    \Phi(\lambda)
    =
    {\rm k}+\lambda{\rm d}	
    +
    \int_{(0,\infty)}(1-e^{-\lambda x})\,\Pi(dx),
    \quad
    \lambda> 0,
    \label{BF}
\end{align}
where ${\rm k, d}\geq0$ and $\Pi$ is a measure on $(0,\infty)$ satisfying $\int_{(0,\infty)}(1\wedge x)\,\Pi(dx)<\infty$. The correspondence between the Bernstein function $\Phi$ and the triplet $({\rm k}, {\rm d},\Pi)$ is one-to-one. A Bernstein function $\Phi\neq 0$ is called \emph{special} if $\lambda/\Phi(\lambda)$ is also a Bernstein function. A Bernstein function is called \emph{complete} if the corresponding measure $\Pi$ admits a completely monotone density $\Pi(dx)/dx$.
A function $\Phi:(0,\infty)\to(0,\infty)$ is a complete Bernstein function if and only if $\lambda/\Phi(\lambda)$ is, as shown in \cite[Proposition 7.1]{SchSonVon}. In particular, any non-zero complete Bernstein function is special.

Let us recall basic facts concerning subordinators. See Bertoin \cite[Chapter III]{Ber}, \cite{Ber99} and Kyprianou \cite[Chapter 5]{Kyp} for more details.
Let $(\Omega, \mathcal{F}, \mathbb{P})$ be a probability space.
Let $\tau=(\tau(s))_{s\geq0}$ be a (killed) \emph{subordinator} with cemetery point $\infty$, that is, a non-negative L\'evy process killed at an independent exponential random time. 
By the L\'evy--Khintchine formula,
the Laplace transform of $\tau(s)$ can be expressed in the form
\begin{align}
   	\mathbb{E}\exp(-\lambda \tau(s))
   	=
   	\exp(-\Phi(\lambda)s),
   	\quad
   	\lambda,s> 0,
   	\label{Laplace-trans}
\end{align}
where it is understood that $\exp(-\infty)=0$, and $\Phi:(0,\infty)\to[0,\infty)$ is a Bernstein function and called the \emph{Laplace exponent} of $\tau$. 
Therefore $\Phi$ is represented as \eqref{BF}, where $\rm k$, $\rm d$ and $\Pi$ are called the \emph{killing rate}, the \emph{drift coefficient} and the \emph{L\'evy measure} of $\tau$, respectively. They are uniquely determined by the distribution of $\tau$. Conversely, for any Bernstein function $\Phi$, there exists a  subordinator $\tau$ satisfying \eqref{Laplace-trans}. A subordinator $\tau$ is called \emph{special} or \emph{complete} if $\Phi$ is special or complete, respectively.

Let 
\begin{align}
	L(t)
	=
	\inf\{s\geq0\::\: \tau(s)\in(t, \infty]\},
	\quad
	t\geq0.
\end{align}
Remark that we regard $\infty$ as the cemetery point of $\tau$, and hence $L(\infty):=\lim_{t\to\infty}L(t)$ is the lifetime of $\tau$.
We are interested in the asymptotic behavior of the undershoot
\begin{align}
    G(t)=
    \tau(L(t)-),
    \quad t\geq0.	
\end{align}
Since $\bP(\tau(L(t)-)<t=\tau(L(t)))=0$ by \cite[Proposition III.2]{Ber}, we have    
\begin{align}
     G(t)=\sup\{\tau(s)\::\: \text{$s\geq0$, $\tau(s)\leq t$}\},
     \quad \text{a.s.}
\end{align}
Note that our convention for the cemetery point differs from that in \cite{Ser25ECP}, ensuring that $G(t)$ is well-defined even if $\{s\geq0\::\: \tau(s)\in(t,\infty)\}=\emptyset$.

From now on we always assume that $\tau$ is not identically zero, or equivalently, $\Phi\neq 0$. 
We denote by $U$ the \emph{potential measure} (or renewal measure) of $\tau$, that is,
\begin{align}
    U(A)
    =
    \bE\int_0^\infty \mathbbm{1}_{\{\tau(s)\in A\}}\,ds,
    \quad
    \text{for any Borel subset $A\subset [0,\infty)$.}	
\end{align}
Then $U$ is a locally finite measure on $[0,\infty)$, and it is diffuse (i.e., $U(\{x\})=0$ for any $x\geq0$) if and only if $\Phi(\lambda)\to \infty$ as $\lambda\to\infty$, or equivalently, $\mathrm{d}>0$ or $\Pi(0,\infty)=\infty$.
The Laplace transform of $U$ is given by
\begin{align}
    \int_{[0,\infty)} e^{-\lambda x}\,U(dx)
    =
    \int_0^\infty \bE\exp(-\lambda \tau(s))\,ds
    =
    \frac{1}{\Phi(\lambda)},
    \quad
    \lambda>0.	
\end{align}
For a measure $\mu$ on $\bR$, we abbreviate $\mu((a,b))$ and $\mu([a,b])$ as $\mu(a,b)$ and $\mu[a,b]$, respectively.

\begin{Prop}\label{Prop:Law-Gt}
Let $t>0$. The following statements hold.
\begin{enumerate}[label=\upshape(\arabic*)]
\item $\bP(G(t)\in dx)
     =
     U(dx)({\rm k}+\Pi(t-x,\infty))$, 
     $0\leq x<t$.
\item If $\mathrm{d}=0$, then
$\bP(G(t)=t)=U(\{t\})({\rm k}+\Pi(0,\infty))$.  Here it is understood that $0\cdot \infty=0$.	
\item Assume $\mathrm{d}>0$. Then $U$ admits a density $u(x)=U(dx)/dx$ that is positive and continuous on $(0,\infty)$. In addition $u(0+)=1/{\rm d}$ and $\bP(G(t)=t)=u(t){\rm d}$.
\end{enumerate}
	
\end{Prop}

For the proof of Proposition \ref{Prop:Law-Gt}, see \cite[Section III.2]{Ber} and \cite[Sections 5.2 and 5.3]{Kyp}.

\begin{Prop}\label{Prop:special-complete}
Assume that a subordinator $\tau$ is not identically zero.
Then the following conditions are equivalent.
\begin{enumerate}[label=\upshape(\roman*)]
\item $\tau$ is special (respectively, complete).
\item The restriction of $U$ over $(0,\infty)$ admits a density $u(x)=U|_{(0,\infty)}(dx)/dx$ that is non-increasing  (respectively, completely monotone) on $(0,\infty)$.	
\end{enumerate}

\end{Prop}

We refer the reader to \cite[ Theorem 11.3 and Remark 11.6]{SchSonVon} for more details of Proposition \ref{Prop:special-complete}.

Let us review the theory of regular variation.
We refer the reader to Bingham, Goldie and Teugels \cite{BinGolTeu} for the details.
A positive measurable function $f$ defined on $(0,R)$ for some $R>0$ is said to be \emph{regularly varying} at $0+$ with index $\rho\in\mathbb{R}$ if $f(\lambda x)/f(x)\to \lambda^\rho$ as $x\to 0+$  for any $\lambda>0$.
Similarly, a positive measurable function $f$ defined on $(R,\infty)$ for some $R\geq0$ is said to be regularly varying at $\infty$ with index $\rho\in\mathbb{R}$ if  $f(\lambda x)/f(x)\to \lambda^\rho$ as $x\to\infty$ for any $\lambda>0$. A regularly varying function with index $0$ is said to be \emph{slowly varying}. Then $f$ is regularly varying at $0+$ (respectively, at $\infty$) with index $\rho$  if and only if $f$ can be represented as $f(x)=x^\rho \ell(x)$ for some function $\ell$ slowly varying at $0+$ (respectively, at $\infty$). 
By Karamata's Tauberian theorem and the monotone density theorem for regularly varying functions, we obtain the following propositions immediately.

\begin{Prop}\label{Prop:RV-long}
Let $0<\alpha<1$. Let $\ell:(0,\infty)\to (0,\infty)$   be slowly varying at $\infty$. Then the following conditions are equivalent:
\begin{enumerate}[label=\upshape(\roman*)]
\item $\Phi(\lambda)
    \sim
    \lambda^\alpha \ell(1/\lambda)$ as $\lambda\to 0+$.

\item ${\rm k}=0$  and $\displaystyle\int_0^x \Pi(y,\infty)\,dy
    \sim
    \frac{x^{1-\alpha}\ell(x)}{\Gamma(2-\alpha)}$ as $x\to\infty$.	

\item ${\rm k}=0$  and $\displaystyle\Pi(x,\infty)
    \sim
    \frac{x^{-\alpha}\ell(x)}{\Gamma(1-\alpha)}$ as $x\to\infty$.

\item $\displaystyle U[0,x]\sim \frac{x^{\alpha}}{\Gamma(1+\alpha)\ell(x)}$ as $x\to \infty$.

\end{enumerate}
Here $\Gamma(z)$ denotes the gamma function.
If, in addition, $\tau$ is special, then the above conditions are also equivalent to
\begin{enumerate}[label=\upshape(\roman*), resume]
\item $\displaystyle u(x)\sim \frac{x^{-1+\alpha}}{\Gamma(\alpha)\ell(x)}$ as $x\to\infty$,
\end{enumerate} 
where $u(x)$ denotes the non-increasing density of $U|_{(0,\infty)}$, whose existence is guaranteed by Proposition {\upshape\ref{Prop:special-complete}}.
\end{Prop}

\begin{Prop}\label{Prop:RV-short}
Let $0<\alpha<1$. Let $\ell:(0,\infty)\to (0,\infty)$   be slowly varying at $0+$. Then the following  conditions are equivalent:
\begin{enumerate}[label=\upshape(\roman*)]
\item\label{Prop:RV-short-1}$\Phi(\lambda)
    \sim
    \lambda^\alpha \ell(1/\lambda)$ as $\lambda\to\infty$.	 
    
\item\label{Prop:RV-short-2} ${\rm d}=0$ and $\displaystyle\int_0^x \Pi(y,\infty)\,dy
    \sim
    \frac{x^{1-\alpha}\ell(x)}{\Gamma(2-\alpha)}$ as $x\to0+$.	

\item\label{Prop:RV-short-3} ${\rm d}=0$ and $\displaystyle\Pi(x,\infty)
    \sim
    \frac{x^{-\alpha}\ell(x)}{\Gamma(1-\alpha)}$  as $x\to0+$.

\item\label{Prop:RV-short-4} $\displaystyle U[0,x]\sim \frac{x^{\alpha}}{\Gamma(1+\alpha)\ell(x)}$  as $x\to0+$.

\end{enumerate}
If, in addition, $\tau$ is special, then the above conditions are also equivalent to
\begin{enumerate}[label=\upshape(\roman*), resume]
\item\label{Prop:RV-short-5} $U(\{0\})=0$ and $\displaystyle u(x)\sim \frac{x^{-1+\alpha}}{\Gamma(\alpha)\ell(x)}$ as $x\to0+$.
\end{enumerate}
\end{Prop}

Let us recall the Dynkin--Lamperti theorem. We refer the reader to \cite[Theorem III.6]{Ber} and \cite[Theorem 5.16]{Kyp} for more details.

\begin{Thm}[Dynkin--Lamperti theorem]\label{Thm:DL}
Let $\alpha\in(0,1)$. The following conditions are equivalent.
\begin{enumerate}[label=\upshape(\roman*)]
	\item $\Phi$ is regularly varying  at $0+$ (respectively, at $\infty$) with index $\alpha$.
\item	 For any $0\leq y \leq 1$, 
\begin{align}
     \mathbb{P}\bigg(\frac{G(t)}{t}\leq y\bigg)
     \to
     \frac{\sin(\pi\alpha)}{\pi}\int_0^y \frac{dx}{x^{1-\alpha}(1-x)^\alpha},
     \label{Thm:long-DL-eq}	
\end{align}
as $t\to\infty$ (respectively, as $t\to0+$).
\end{enumerate}
\end{Thm}
Note that the right-hand side of \eqref{Thm:long-DL-eq} is the distribution function of the Beta distribution $B(\alpha,1-\alpha)$, which is the usual arcsine distribution in the case $\alpha=1/2$.



\section{Large Deviation Estimates Related to Dynkin--Lamperti Theorem}\label{Sec:Main}

We are interested in asymptotic behavior of the left-hand side of \eqref{Thm:long-DL-eq}  as $t\to\infty$ (or $t\to 0+$) and $y\to0+$.
For this purpose,
we need mild versions of Propositions \ref{Prop:RV-long} and \ref{Prop:RV-short} as auxiliary results. 

\begin{Prop}\label{Prop:mild-RV-long}
Let $\alpha\in(0,1)$. Let $\ell:(0,\infty)\to(0,\infty)$ be slowly varying at $\infty$.
Then the following conditions are equivalent.
\begin{enumerate}[label=\upshape(\roman*)]
\setlength{\itemsep}{0.5mm}
 \item \label{Prop:mild-RV-long-item-1}	$\Phi(\lambda)\asymp \lambda^\alpha \ell(1/\lambda)$ as $\lambda\to0+$.
 \item\label{Prop:mild-RV-long-item-2} ${\rm k}=0$ and $\int_0^x \Pi(y,\infty)\,dy\asymp x^{1-\alpha}\ell(x)$ as $x\to\infty$.
 \item\label{Prop:mild-RV-long-item-3} ${\rm k}=0$ and $\Pi(x,\infty)\asymp x^{-\alpha}\ell(x)$ as $x\to\infty$. 
 \item\label{Prop:mild-RV-long-item-4} $U[0,x]\asymp x^{\alpha}/\ell(x)$ as $x\to\infty$.
\end{enumerate}
If, in addition, $\tau$ is special, then the above conditions are also equivalent to
\begin{enumerate}[resume, label=\upshape(\roman*)]
 \item\label{Prop:mild-RV-long-item-5} $u(x)\asymp x^{-1+\alpha}/\ell(x)$ as $x\to\infty$.	
\end{enumerate}
\end{Prop}

\begin{Prop}\label{Prop:mild-RV-short}
Let $\alpha\in(0,1)$. Let $\ell:(0,\infty)\to(0,\infty)$ be slowly varying at $0+$.
Then the following conditions are equivalent.
\begin{enumerate}[label=\upshape(\roman*)]
\setlength{\itemsep}{0.5mm}
 \item\label{Prop:mild-RV-short-item-1} 	$\Phi(\lambda)\asymp \lambda^\alpha \ell(1/\lambda)$ as $\lambda\to\infty$.
 \item\label{Prop:mild-RV-short-item-2}  ${\rm d}=0$ and $\int_0^x \Pi(y,\infty)\,dy\asymp x^{1-\alpha}\ell(x)$ as $x\to0+$.
 \item\label{Prop:mild-RV-short-item-3}  ${\rm d}=0$ and $\Pi(x,\infty)\asymp x^{-\alpha}\ell(x)$ as $x\to0+$. 
 \item\label{Prop:mild-RV-short-item-4}  $U[0,x]\asymp x^{\alpha}/\ell(x)$ as $x\to0+$.
\end{enumerate}
If, in addition, $\tau$ is special, then the above conditions are also equivalent to
\begin{enumerate}[resume, label=\upshape(\roman*)]
 \item\label{Prop:mild-RV-short-item-5}  $U(\{0\})=0$ and $u(x)\asymp x^{-1+\alpha}/\ell(x)$ as $x\to0+$.	
\end{enumerate}
\end{Prop}

The proofs of Propositions \ref{Prop:mild-RV-long}
and \ref{Prop:mild-RV-short} 
 will be given in Section \ref{Sec:Tauberian-type}. 
We remark that the equivalence between (i), (ii) and (iv) follows from a Tauberian theorem for dominated variations (see
de Haan and Stadtm\"uller \cite[Theorem 1]{deHSta} and  
 Bingham, Goldie, and Teugels \cite[Theorem 2.10.2]{BinGolTeu}). See also Bertoin \cite[Proposition III.1]{Ber}. Nevertheless, the corresponding monotone density theorem for dominated variations in \cite[Theorem 2.10.2]{BinGolTeu} does not yield [(ii) $\Rightarrow$ (iii)] nor [(iv) $\Rightarrow$ (v)]. Instead, we follow the methods of  Z\"ahle \cite[Theorem 7]{Zah} and Kim, Song, and Vondra{\v c}ek \cite[Section 2.2]{KimSongVon}, which imply Proposition \ref{Prop:mild-RV-short} immediately. Nevertheless, their approaches seem not to work for the proof of Proposition \ref{Prop:mild-RV-long}. So we will improve their methods to  prove Propositions \ref{Prop:mild-RV-long} and \ref{Prop:mild-RV-short} in almost the same way.

By using the above propositions, we establish large deviation estimates related to the Dynkin--Lamperti theorem as follows.

\begin{Thm}[Long-range large deviation estimates]\label{Thm:LDE-long}
Let $\alpha\in(0,1)$ and $\hat{\alpha}=1-\alpha$. Let $\ell:(0,\infty)\to(0,\infty)$ be slowly varying at $\infty$, and suppose $\varepsilon:(0,\infty)\to(0,1)$ satisfies $\varepsilon(t)\to 0$ and $\varepsilon(t)t\to\infty$ as $t\to\infty$.
 Then the following statements hold.

\begin{enumerate}[label=\upshape(\arabic*)]
\item\label{Thm:LDE-long-item-1} Assume $\Phi(\lambda)\sim \lambda^\alpha \ell(1/\lambda)$ as $\lambda\to0+$. Then
\begin{align}\label{Thm:LDE-long-1}
  	 \mathbb{P}\bigg(\frac{G(t)}{t}\leq \varepsilon(t)\bigg)
     \sim
     \frac{\sin(\pi\alpha)}{\pi\alpha}
     \frac{\ell(t)}{\ell(\varepsilon(t)t)}
    \varepsilon(t)^\alpha,
\end{align}
as $t\to\infty$. Suppose additionally that $\tau$ is special. Then
\begin{align}\label{Thm:LDE-long-2}
	 \mathbb{P}\bigg(\frac{G(t)}{t}\geq 1- \varepsilon(t)\bigg)
     \sim
     \frac{\sin(\pi\hat{\alpha})}{\pi\hat{\alpha}}
     \frac{\ell(\varepsilon(t)t)}{\ell(t)}
     \varepsilon(t)^{\hat{\alpha}},
\end{align}
as $t\to\infty$.

\item\label{Thm:LDE-long-item-2} Assume $\Phi(\lambda)\asymp \lambda^\alpha \ell(1/\lambda)$ as $\lambda\to0+$. Then
\begin{align}\label{Thm:LDE-long-3}
  	 \mathbb{P}\bigg(\frac{G(t)}{t}\leq \varepsilon(t)\bigg)
     \asymp
     \frac{\ell(t)}{\ell(\varepsilon(t)t)}
     \varepsilon(t)^\alpha,
\end{align}
as $t\to\infty$.
Suppose additionally that $\tau$ is special. Then
\begin{align}\label{Thm:LDE-long-4}
	 \mathbb{P}\bigg(\frac{G(t)}{t}\geq 1- \varepsilon(t)\bigg)
     \asymp
     \frac{\ell(\varepsilon(t)t)}{\ell(t)}
     \varepsilon(t)^{\hat{\alpha}},
\end{align}
as $t\to\infty$.

\end{enumerate}
\end{Thm}

The following theorem is a short-range version of Theorem \ref{Thm:LDE-long}.

\begin{Thm}[Short-range large deviation estimates]\label{Thm:LDE-short}
Let $\alpha\in(0,1)$ and $\hat{\alpha}=1-\alpha$. Let $\ell:(0,\infty)\to(0,\infty)$ be slowly varying at $0+$, and suppose $\varepsilon:(0,\infty)\to(0,1)$ satisfies $\varepsilon(t)\to 0$ as $t\to0+$.
Then the following statements hold.
\begin{enumerate}[label=\upshape(\arabic*)]
\item Assume $\Phi(\lambda)\sim \lambda^\alpha \ell(1/\lambda)$ as $\lambda\to\infty$, then \eqref{Thm:LDE-long-1} holds as $t\to0+$. Suppose additionally that $\tau$ is special. Then
\eqref{Thm:LDE-long-2} holds as $t\to0+$.

\item Assume $\Phi(\lambda)\asymp \lambda^\alpha \ell(1/\lambda)$ as $\lambda\to\infty$, then \eqref{Thm:LDE-long-3} holds as $t\to0+$. Suppose additionally that $\tau$ is special. Then
\eqref{Thm:LDE-long-4} holds as $t\to0+$.

\end{enumerate}
\end{Thm}

The proofs of Theorems \ref{Thm:LDE-long} and \ref{Thm:LDE-short} will be given in 
Section \ref{sec:Proof-Main}.

\begin{Rem}[Comparison with the Dynkin--Lamperti theorem]
Note that
\begin{align}
   	 \frac{\sin(\pi\alpha)}{\pi}\int_0^{y} \frac{dx}{x^{1-\alpha}(1-x)^\alpha}
   	\sim
   	 \frac{\sin(\pi\alpha)}{\pi}\int_0^{y} \frac{dx}{x^{1-\alpha}}
   	 =
   	 \frac{\sin(\pi\alpha)}{\pi\alpha}y^\alpha,
\quad
\text{as $y\to0+$},
\end{align}
and
\begin{align}
   	 \frac{\sin(\pi\alpha)}{\pi}\int_{1-y}^{1} \frac{dx}{x^{1-\alpha}(1-x)^\alpha}
   	\sim
   	 \frac{\sin(\pi\alpha)}{\pi}\int_{1-y}^{1} \frac{dx}{(1-x)^{\alpha}}
   	 =
   	 \frac{\sin(\pi\hat{\alpha})}{\pi\hat{\alpha}}y^{\hat{\alpha}},
\quad
\text{as $y\to0+$}.
\end{align}
Therefore, Theorems \ref{Thm:LDE-long} and \ref{Thm:LDE-short} are consistent with Theorem \ref{Thm:DL} if we formally exchange the order in which the limits are taken and regard $\ell(t)/\ell(\varepsilon(t)t)$ as constant.

\end{Rem}


\section{Tauberian-type Estimates for Laplace Transforms}
\label{Sec:Tauberian-type}

For the proofs of Propositions \ref{Prop:mild-RV-long} and \ref{Prop:mild-RV-short}, 
we will need the following Tauberian-type estimates for Laplace transforms of monotone functions, without assuming regular variation conditions. 

\begin{Lem}\label{Lem:Lap-inv-increasing}
Let $v:(0,\infty)\to[0,\infty)$ be a non-decreasing function with
\begin{align}
      \cL v(\lambda)
      =
      \int_0^\infty e^{-\lambda u}v(u)\,du\in(0,\infty),
      \quad \lambda>0.
      \label{Laplace-w-v}
\end{align}
Then the following statements hold.
\begin{enumerate}[label=\upshape(\arabic*)]
\item $\displaystyle \frac{e^{y}}{x}\cL v\bigg(\frac{1}{x}\bigg)\geq v(xy)$ for any $x,y>0$.
\item  Assume that $C_0=\inf_{\lambda>0} (\cL v(2\lambda)/\cL v(\lambda))>0$.
Then there exists $C>0$ such that
\begin{align}
	\frac{C}{x} \cL v\bigg(\frac{1}{x}\bigg)
	\leq
	v(x),
	\quad \text{for all $x>0$.}
\end{align}
\end{enumerate}
\end{Lem}

\begin{proof}
We follow the proof of Z\"ahle \cite[Theorem 7]{Zah}.

(1) Since $v$ is non-decreasing, we have
\begin{align}
  \frac{1}{x}\cL v\bigg(\frac{1}{x}\bigg)
  &=
  \frac{1}{x}
  \int_0^\infty e^{-u/x}v(u)\,du
  =
  \int_0^\infty e^{-u}v(xu)\,du
  \\
  &\geq
  \bigg(\int_y^\infty e^{-u}\,du\bigg)v(xy)=e^{-y} v(xy),	
\end{align}
as desired.

(2) Let $R>0$. By using (1) and the assumption, we have
\begin{align}
   \frac{R}{x}\cL v\bigg(\frac{R}{x}\bigg)
   &=
   \int_0^R +\int_R^\infty e^{-y}v(xy/R)\,dy
   \\
   &\leq
   \bigg(\int_0^R e^{-y}\,dy\bigg) v(x)
   +
   \int_R^\infty e^{-y}
   \frac{R e^{y/2}}{2x}\cL v\bigg(\frac{R}{2x}\bigg)\,dy
   \\
   &=
   (1-e^{-R})v(x)
   +
   \frac{R e^{-R/2}}{x}\cL v\bigg(\frac{R}{2x}\bigg)
   \\
   &\leq
   v(x)
   +
   \frac{R C_0 e^{-R/2}}{x}\cL v\bigg(\frac{R}{x}\bigg),
   \quad
   x>0.
\end{align}
Let us take $R=2^{n}$  for some $n\in\bN$ large enough so that $C_0 e^{-R/2}<1/2$. Then
\begin{align}
	v(x)\geq (1-C_0 e^{-R/2})\frac{R}{x}\cL v\bigg(\frac{R}{x}\bigg)
	\geq \frac{R}{2x}\cL v\bigg(\frac{R}{x}\bigg)
	\geq \frac{RC_0^{n}}{2x}\cL v\bigg(\frac{1}{x}\bigg),
	\quad
	x>0,
\end{align}
as desired. 
\end{proof}

\begin{Lem}\label{Lem:Lap-inv-decreasing}
Let $v:(0,\infty)\to[0,\infty)$ be a non-increasing function with \eqref{Laplace-w-v}.
Set
\begin{align}
    V(x)=\int_0^x v(y)\,dy,\quad x>0.	
\end{align}
Then the following statements hold.
\begin{enumerate}[label=\upshape(\arabic*)]
\item\label{Lem:Lap-inv-decreasing-1} $\displaystyle \frac{1}{1-e^{-1}}\frac{1}{x}\cL v\bigg(\frac{1}{x}\bigg)\geq v(x)$ for any $x>0$.
\item\label{Lem:Lap-inv-decreasing-2} 
$\displaystyle e\cL v\bigg(\frac{1}{x}\bigg)\geq V(x)$ for any $x>0$.

\item\label{Lem:Lap-inv-decreasing-3} 
Assume that there exists $R_0>0$ such that
\begin{align}\label{Asmp:ratio-bdd-power-v3}
  f(\lambda)=\sup_{0<s\leq R_0}\frac{\cL v(s)}{\cL v(\lambda s)}\to 0,
  \quad
  \text{as $\lambda\to 0+$}.	
\end{align}
Then there exists $C,R>0$ such that
\begin{align}
  \frac{C}{x}\cL v\bigg(\frac{1}{x}\bigg)\leq v(x),
  \quad
  \text{for any $x>R$.}	
\end{align}

\item\label{Lem:Lap-inv-decreasing-4} 
Assume that there exists $R_0>0$ such that
\begin{align}\label{Asmp:ratio-bdd-power}
  g(\lambda)=\sup_{s\geq R_0}\frac{\cL v(\lambda s)}{\cL v(s)}\to 0,
  \quad
  \text{as $\lambda\to \infty$}.	
\end{align}
Then there exists $C,R>0$ such that
\begin{align}
  \frac{C}{x}\cL v\bigg(\frac{1}{x}\bigg)\leq v(x),
  \quad
  \text{for any $0<x<R$.}	
\end{align}

\end{enumerate}
\end{Lem}

We remark that  a mild version of \ref{Lem:Lap-inv-decreasing-4} can also be found in Z\"ahle \cite[Theorem 7]{Zah} and Kim, Song, and Vondra{\v c}ek \cite[Proposition 2.5]{KimSongVon}. Nevertheless, their approaches seem not to work for the proof of  \ref{Lem:Lap-inv-decreasing-3}. In the following, we will provide a unified method for proving  \ref{Lem:Lap-inv-decreasing-3} and \ref{Lem:Lap-inv-decreasing-4} simultaneously.

\begin{proof}
\ref{Lem:Lap-inv-decreasing-1} Since $v$ is non-increasing, we have
\begin{align}
  \frac{1}{x}\cL v\bigg(\frac{1}{x}\bigg)
  &=
  \frac{1}{x}\int_0^\infty e^{- y/x}v(y)\,dy
  =\int_0^\infty e^{-y}v(xy)\,dy
  \\
  &\geq 
  \bigg(\int_0^1 e^{-y}\,dy\bigg)v(x)=(1-e^{-1})v(x),	
\end{align}
as desired.	

\ref{Lem:Lap-inv-decreasing-2} By the definition, $V(x)$ is non-decreasing. We use Fubini's theorem to see
\begin{align}
   \cL v(\lambda)
    &=
    \int_0^\infty e^{-\lambda x}v(x)\,dx
    =
    \int_0^\infty \bigg(\int_x^\infty \lambda e^{-\lambda y}\,dy\bigg) v(x)\,dx
    \\
    &=
    \lambda\int_0^\infty e^{-\lambda y}\bigg(\int_0^y v(x)\,dx\bigg)\,dy
    =
    \lambda \int_0^\infty e^{-\lambda y}V(y)\,dy,	
\end{align}
in other words,
\begin{align}
    \frac{\cL v(\lambda)}{\lambda}=\int_0^\infty e^{-\lambda x}V(x)\,dx.	
\end{align}
Therefore Lemma \ref{Lem:Lap-inv-increasing} implies
\begin{align}
  e\cL v\bigg(\frac{1}{x}\bigg)=\frac{e}{x}\frac{\cL v(1/x)}{1/x}
  \geq
  V(x)
  \quad
  \text{for any $x>0$,}	
\end{align}
as desired.

\ref{Lem:Lap-inv-decreasing-3}
Let $\varepsilon\in(0,1]$. We use the fact that $v(x)$ is non-increasing and \ref{Lem:Lap-inv-decreasing-2} to get
\begin{align}
	\cL v\bigg(\frac{\varepsilon}{x}\bigg)
	&=
    \int_0^x +\int_{x}^\infty e^{-\varepsilon y/x}v(y)\,dy
	\leq
	V(x)+\frac{x}{\varepsilon}v(x)
	\leq
	e \cL v\bigg(\frac{1}{x}\bigg)+\frac{x}{\varepsilon}v(x).
	\label{lower-bdd-proof-v2-1} 
\end{align}
By the definition of $f$, we have
\begin{align}
   \frac{1}{f(\varepsilon)}\cL v\bigg(\frac{1}{x}\bigg)
   \leq
   \cL v\bigg(\frac{\varepsilon}{x}\bigg),
   \quad
   x>\frac{1}{R_0}.	
   \label{lower-bdd-proof-v2-2} 
\end{align}
Combining \eqref{lower-bdd-proof-v2-1} with \eqref{lower-bdd-proof-v2-2}, we get
\begin{align}
  \frac{1}{f(\varepsilon)}\cL v\bigg(\frac{1}{x}\bigg)
   \leq
  e \cL v\bigg(\frac{1}{x}\bigg)+\frac{x}{\varepsilon}v(x),
  \quad
  x>\frac{1}{R_0}.	
\end{align}
Take $\varepsilon\in(0,1]$ small enough so that $1/f(\varepsilon)\geq e+1$. Then
\begin{align}
    \frac{\varepsilon}{x}\cL v\bigg(\frac{1}{x}\bigg)
    \leq 
    v(x),
    \quad x>\frac{1}{R_0}, 	
\end{align}
as desired.

\ref{Lem:Lap-inv-decreasing-4}
We can obtain the desired result by replacing $f(\varepsilon)$ and $x>1/R_0$ with $g(1/\varepsilon)$ and $0<x<\varepsilon/R_0$, respectively, in the proof of \ref{Lem:Lap-inv-decreasing-3}.
\end{proof}

Let us apply the above lemmas to  subordinators.

\begin{Cor}\label{Cor:bdd}
Let $\tau$ be a subordinator with Laplace exponent $\Phi\neq0$. Then
 \begin{align}
    \Phi\bigg(\frac{1}{x}\bigg)x
  	\asymp
  	{\rm k}x+{\rm d}+\int_0^x \Pi(y,\infty)\,dy,
  	\quad
  	\text{for all $x>0$},
  	\label{Cor:bdd-1}
\end{align}
and
\begin{align}
  	\frac{1}{\Phi(1/x)}
  	\asymp
  	U[0,x],
  	\quad
  	\text{for all $x>0$}.
  	\label{Cor:bdd-2}
 \end{align}
\end{Cor}

This kind of estimates can also be found in Bertoin \cite[Proposition III.1]{Ber}. Nevertheless we prove it as a corollary of Lemma \ref{Lem:Lap-inv-increasing} in order to clarify the importance of the fact that the right hand sides of \eqref{Cor:bdd-1} and \eqref{Cor:bdd-2} are non-decreasing functions. 

\begin{proof}

Note that
\begin{align}
   \frac{\Phi(\lambda)}{\lambda^2}
   &=
   \frac{1}{\lambda^2}
   \bigg({\rm k}+\lambda{\rm d}
   +\lambda\int_0^\infty e^{-\lambda x}\Pi(x,\infty)\,dx 
   \bigg)
   \\
   &=
   \int_0^\infty e^{-\lambda x}\bigg({\rm k}x+{\rm d}+\int_0^x \Pi(y,\infty)\,dy\bigg)\,dx,
   \quad
   \lambda>0.	
\end{align}
In addition, $\Phi$ is non-decreasing and hence
\begin{align}
    	\frac{\Phi(2\lambda)/(2\lambda)^2}{\Phi(\lambda)/\lambda^2}\geq \frac{1}{4},
    	\quad
    	\lambda>0.
\end{align}
So we use Lemma \ref{Lem:Lap-inv-increasing} with $\cL v(\lambda)=\Phi(\lambda)/\lambda^2$ and $v(x)={\rm k}x+{\rm d}+\int_0^x \Pi(y,\infty)\,dy$ to obtain 
\eqref{Cor:bdd-1}. 	

Let us prove \eqref{Cor:bdd-2}.  Note that
\begin{align}
    \frac{1}{\lambda \Phi(\lambda)}
    =
    \frac{1}{\lambda} \int_{[0,\infty)} e^{-\lambda x}\,U(dx)
    =
    \int_0^\infty e^{-\lambda x} U[0,x]\,dx,
    \quad
    \lambda>0. 	
\end{align}
Since $\Phi$ is a Bernstein function, it is concave and non-negative with $\Phi(0+)={\rm k}$, which implies $(\Phi(x)+\Phi(y))/2\leq \Phi((x+y)/2)$ for all $x,y>0$. In particular, 
\begin{align}
    \frac{\Phi(2\lambda)}{2}\leq\frac{\Phi(2\lambda)+\Phi(0+)}{2}\leq \Phi(\lambda),
    \quad
    \lambda>0,	
\end{align}
and hence $(\lambda\Phi(\lambda))/(2\lambda\Phi(2\lambda))\geq 1/4$.
Applying Lemma \ref{Lem:Lap-inv-increasing} to $\cL v(\lambda)=1/(\lambda\Phi(\lambda))$ and $v(x)=U[0,x]$, we obtain \eqref{Cor:bdd-2}.
\end{proof}

We are ready to proceed with the proof of Proposition \ref{Prop:mild-RV-long}.

\begin{proof}[Proof of Proposition \upshape\ref{Prop:mild-RV-long}]
The equivalence between \ref{Prop:mild-RV-long-item-1}, \ref{Prop:mild-RV-long-item-2} and \ref{Prop:mild-RV-long-item-4} follows immediately from Corollary \ref{Cor:bdd}.

Proof of \ref{Prop:mild-RV-long-item-1} $\Rightarrow$ \ref{Prop:mild-RV-long-item-3}. Note that
\begin{align}
\frac{\Phi(\lambda)}{\lambda}
   -{\rm d}
   =
   \int_0^\infty e^{-\lambda x}\Big({\rm k}+\Pi(x,\infty)\Big)\,dx,
   \quad
   \lambda>0.	
\end{align}
Assume $\Phi(\lambda)\asymp \lambda^\alpha \ell(1/\lambda)$ as $\lambda\to 0+$. Let $\varepsilon\in(0, 1-\alpha)$. By Potter's theorem \cite[Theorem 1.5.6]{BinGolTeu}, we can take  $R_0>0$ so that
\begin{align}
    \frac{\ell(1/s)}{\ell(1/(\lambda s))}\leq 2 \lambda^{-\varepsilon},
    \quad
    \text{for any $0<\lambda\leq 1$ and $0<s<R_0$.}	
\end{align}
Then $\Phi(s)/s\to \infty$ as $s\to 0+$ and
\begin{align}
   \frac{(\Phi(s)/s)-{\rm d}}{(\Phi(\lambda s)/(\lambda s))-{\rm d}}
   \sim
   \frac{\Phi(s)/s}{\Phi(\lambda s)/(\lambda s)}
   \asymp
  \lambda^{1-\alpha}\frac{\ell(1/s)}{\ell(1/(\lambda s))} 
  \leq 2\lambda^{1-\alpha-\varepsilon},	
\end{align}
uniformly in $\lambda\in(0,1]$ as $s\to0+$. Therefore we use Lemma \ref{Lem:Lap-inv-decreasing} \ref{Lem:Lap-inv-decreasing-1} and \ref{Lem:Lap-inv-decreasing-3} with $\cL v(\lambda)=(\Phi(\lambda)/\lambda)-{\rm d}$ and $v(x)={\rm k}+\Pi(x,\infty)$ to obtain ${\rm k}=0$ and $\Pi(x,\infty)\asymp x^{-\alpha}\ell(x)$ as $x\to\infty$.

Proof of \ref{Prop:mild-RV-long-item-3} $\Rightarrow$ \ref{Prop:mild-RV-long-item-2}. Assume $\Pi(x,\infty)\asymp x^{-\alpha}\ell(x)$ as $x\to\infty$. By \cite[Corollary 1.4.2]{BinGolTeu}, the slowly varying function $\ell$ is locally bounded and locally integrable on $[R,\infty)$ for some $R>0$. By Karamata's theorem \cite[Theorem 1.5.11]{BinGolTeu}, we have
\begin{align}
    \int_R^x y^{-\alpha}\ell(y)\,dy\sim \frac{x^{1-\alpha}\ell(x)}{1-\alpha},
    \quad
    \text{as $x\to\infty$},	
\end{align}
which immediately implies $\int_0^x \Pi(y,\infty)\,dy\asymp x^{1-\alpha}\ell(x)$ as $x\to\infty$, since $\int_0^R \Pi(y,\infty)\,dy<\infty$.

Proof of \ref{Prop:mild-RV-long-item-1} $\Rightarrow$ \ref{Prop:mild-RV-long-item-5} in the case $\tau$ is special. Recall that $u(x)$ is a density of $U|_{(0,\infty)}$ and it is non-increasing on $(0,\infty)$.
Note that
\begin{align}
    \frac{1}{\Phi(\lambda)}
    -
    U(\{0\})
    =
    \int_0^\infty e^{-\lambda x}u(x)\,dx,
    \quad
    \lambda>0.	
\end{align}
Assume $\Phi(\lambda)\asymp \lambda^\alpha \ell(1/\lambda)$ as $\lambda\to 0+$. Let $\varepsilon\in(0,\alpha)$. Note that $1/\ell$ is slowly varying at $\infty$. By Potter's theorem, we can take $R_0>0$ so that
\begin{align}
    \frac{\ell(1/s)}{\ell(1/(\lambda s))}\geq \frac{ \lambda^{\varepsilon}}{2},
    \quad
    \text{for any $0<\lambda\leq 1$ and $0<s<R_0$.}	
\end{align}
Then $1/\Phi(s)\to\infty$ as $s\to 0+$
and
\begin{align}
   \frac{(1/\Phi(s))-U(\{0\})}{(1/\Phi(\lambda s))-U(\{0\})}
   \sim
   \frac{\Phi(\lambda s)}{\Phi(s)}
   \asymp
   \lambda^\alpha\frac{\ell(1/(\lambda s))}{\ell(1/s)}
   \leq
   2\lambda^{\alpha-\varepsilon},	
\end{align}
uniformly in $\lambda\in(0,1]$ as $s\to 0+$. Therefore we use Lemma \ref{Lem:Lap-inv-decreasing} \ref{Lem:Lap-inv-decreasing-1} and \ref{Lem:Lap-inv-decreasing-3} with $\cL v(\lambda)=(1/\Phi(\lambda))-U(\{0\})$ and $v(x)=u(x)$ to obtain $u(x)\asymp x^{-1+\alpha}/\ell(x)$ as $x\to\infty$.

Proof of \ref{Prop:mild-RV-long-item-5} $\Rightarrow$ \ref{Prop:mild-RV-long-item-4} in the case $\tau$ is special. 
Assume $u(x)\asymp x^{-1+\alpha}/\ell(x)$ as $x\to\infty$. By \cite[Corollary 1.4.2]{BinGolTeu}, we can take $R>0$ large enough so that $1/\ell$ is locally bounded and locally integrable on $[R,\infty)$. By Karamata's theorem, we have
\begin{align}
    \int_R^x \frac{y^{-1+\alpha}}{\ell(y)}\,dy
    \sim
    \frac{x^{\alpha}}{\alpha \ell(x)},
    \quad
    \text{as $x\to\infty$,}	
\end{align}
which immediately implies $U[0,x]\asymp x^\alpha/\ell(x)$ as $x\to\infty$, since $U[0,R]=U(\{0\})+\int_0^R u(y)\,dy<\infty$.
\end{proof}

The proof  of Proposition \ref{Prop:mild-RV-short} is almost the same as that of Proposition \ref{Prop:mild-RV-long}, as we shall see in the following.

\begin{proof}[Proof of Proposition \upshape\ref{Prop:mild-RV-short}]
The equivalence between \ref{Prop:mild-RV-long-item-1}, \ref{Prop:mild-RV-long-item-2} and \ref{Prop:mild-RV-long-item-4} follows immediately from Corollary \ref{Cor:bdd}.

Proof of  \ref{Prop:mild-RV-short-item-1} $\Rightarrow$ \ref{Prop:mild-RV-short-item-3}.
Assume $\Phi(\lambda)\asymp \lambda^\alpha \ell(1/\lambda)$ as $\lambda\to \infty$.  
In this situation the drift coefficient ${\rm d}$ must be zero, and hence
\begin{align}
\frac{\Phi(\lambda)}{\lambda}
   =
   \int_0^\infty e^{-\lambda x}\Big({\rm k}+\Pi(x,\infty)\Big)\,dx,
   \quad
   \lambda>0.	
\end{align}
Let $\varepsilon\in(0,1-\alpha)$. 
Note that $\ell(1/x)$ is slowly varying at $x=\infty$.
By Potter's theorem, we can take $R_0>0$ so that
\begin{align}
    \frac{\ell(1/(\lambda s))}{\ell(1/s)}
    \leq
    2 \lambda^{\varepsilon},
    \quad
    \text{for any $\lambda\geq1$ and $s>R_0$.}
\end{align}
Then
\begin{align}
     \frac{\Phi(\lambda s)/(\lambda s)}{\Phi(s)/s}
   \asymp
  \lambda^{-1+\alpha}\frac{\ell(1/(\lambda s))}{\ell(1/s)} 
  \leq 2\lambda^{-1+\alpha+\varepsilon},		
\end{align}
uniformly in $\lambda\in[1,\infty)$ as $s\to\infty$. Therefore we use Lemma \ref{Lem:Lap-inv-decreasing} \ref{Lem:Lap-inv-decreasing-1} and \ref{Lem:Lap-inv-decreasing-4} with $\cL v(\lambda)=\Phi(\lambda)/\lambda$ and $v(x)={\rm k}+\Pi(x,\infty)$ to obtain $\Pi(x,\infty)\sim v(x)\asymp x^{-\alpha}\ell(x)$ as $x\to0+$, as desired.

Proof of \ref{Prop:mild-RV-short-item-3} $\Rightarrow$ \ref{Prop:mild-RV-short-item-2}. Assume $\Pi(x,\infty)\asymp x^{-\alpha}\ell(x)$ as $x\to0+$. Note that $\ell(1/x)$ is slowly varying at $x=\infty$. By Karamata's theorem, we have
\begin{align}
    \int_0^x \Pi(y,\infty)\,dy
    \asymp\int_0^x y^{-\alpha}\ell(y)\,dy
    =
    \int_{1/x}^\infty y^{-2+\alpha}\ell(1/y)\,dy\sim \frac{x^{1-\alpha}\ell(x)}{1-\alpha},
    \quad
    \text{as $x\to0+$}.	
\end{align}

Proof of \ref{Prop:mild-RV-short-item-1} $\Rightarrow$ \ref{Prop:mild-RV-short-item-5} in the case $\tau$ is special. Assume $\Phi(\lambda)\asymp \lambda^\alpha \ell(1/\lambda)$ as $\lambda\to \infty$. In this case $U(\{0\})=0$ and hence
\begin{align}
    \frac{1}{\Phi(\lambda)}
    =
    \int_0^\infty e^{-\lambda x}u(x)\,dx,
    \quad
    \lambda>0.	
\end{align}
Let $\varepsilon\in(0,\alpha)$.  By Potter's theorem, we can take $R_0>0$ so that
\begin{align}
    \frac{\ell(1/(\lambda s))}{\ell(1/s)}\geq \frac{ \lambda^{-\varepsilon}}{2},
    \quad
    \text{for any $\lambda\geq 1$ and $s>R_0$.}	
\end{align}
Hence
\begin{align}
   \frac{1/\Phi(\lambda s)}{1/\Phi(s)}
   \asymp
   \lambda^{-\alpha}\frac{\ell(1/s)}{\ell(1/(\lambda s))}
   \leq
   2\lambda^{-\alpha+\varepsilon},	
\end{align}
uniformly in $\lambda\in[1,\infty)$ as $s\to \infty$. Therefore we use Lemma \ref{Lem:Lap-inv-decreasing} \ref{Lem:Lap-inv-decreasing-1} and \ref{Lem:Lap-inv-decreasing-4} with $w(\lambda)=1/\Phi(\lambda)$ and $v(x)=u(x)$ to obtain $u(x)\asymp x^{-1+\alpha}/\ell(x)$ as $x\to0+$.

Proof of \ref{Prop:mild-RV-short-item-5} $\Rightarrow$ \ref{Prop:mild-RV-short-item-4} in the case $\tau$ is special. 
Assume $U(\{0\})=0$ and $u(x)\asymp x^{-1+\alpha}/\ell(x)$ as $x\to0+$. By Karamata's theorem, we have
\begin{align}
    U[0,x]
    =
    \int_0^x u(y)\,dy
    \asymp
    \int_0^x \frac{y^{-1+\alpha}}{\ell(y)}\,dy
    \sim
    \frac{x^{\alpha}}{\alpha \ell(x)},
    \quad
    \text{as $x\to0+$.}	
\end{align}
\end{proof}

\section{Proofs of Large Deviation Estimates}\label{sec:Proof-Main}

In this section we will prove Theorems \ref{Thm:LDE-long} and \ref{Thm:LDE-short}.

\begin{Lem}\label{Lem:Gt-estimate}
Let $\tau$ be a subordinator with Laplace exponent $\Phi\neq0$. Then, for any $c\in(0,1)$ and $t>0$,
\begin{align}
   U[0,ct]\Big({\rm k}+\Pi(t,\infty)\Big)
   \leq
   \bP\bigg(\frac{G(t)}{t}\leq c\bigg)
   \leq 
   U[0,ct]\Big({\rm k}+\Pi(t-ct,\infty)\Big).
   \label{Lem:Gt-estimate-1}
\end{align}
If, in addition, $\tau$ is special, then
\begin{align}
	u(t)\bigg({\rm d}+\int_0^{ct}\Big({\rm k}+\Pi(x,\infty)\Big)\,dx\bigg)
	&\leq
	\bP\bigg(\frac{G(t)}{t}\geq 1-c\bigg)
	\\
	&\leq
	u(t-ct)\bigg({\rm d}+\int_0^{ct}\Big({\rm k}+\Pi(x,\infty)\Big)\,dx\bigg).
	\label{Lem:Gt-estimate-2}
\end{align}
\end{Lem}

\begin{proof}
By Proposition \ref{Prop:Law-Gt}, 
\begin{align}
     \bP\bigg(\frac{G(t)}{t}\leq c\bigg)
     =
     \int_{[0,ct]}
     U(dx)
     \Big({\rm k}+\Pi(t-x,\infty)\Big).
     \label{Lem:Gt-estimate-3}	
\end{align}
Since the map $[0,ct]\ni x\mapsto \Pi(t-x,\infty)\in[0,\infty)$ is non-decreasing, we obtain \eqref{Lem:Gt-estimate-1}.	

Further assume that $\tau$ is special. We use Propositions \ref{Prop:Law-Gt} and \ref{Prop:special-complete} to see
\begin{align}
  \bP\bigg(\frac{G(t)}{t}\geq 1-c\bigg)
  &=
  u(t){\rm d}
  +
  \int_{t-ct}^t
  u(x)\Big({\rm k}+\Pi(t-x, \infty)\Big)\,dx
  \\
  &=
  u(t){\rm d}
  +\int_{0}^{ct}
  u(t-x)
  \Big({\rm k}+\Pi(x,\infty)\Big)\,dx
  .
  \label{Lem:Gt-estimate-4}
\end{align}
Note that \eqref{Lem:Gt-estimate-4} is valid even if $\rm d=0$, since $U(\{t\})=0$ in this situation. 
Since the map $[0,ct]\ni x\mapsto u(t-x)\in[0,\infty)$ is non-decreasing, we obtain \eqref{Lem:Gt-estimate-2} as desired.
\end{proof}

\begin{proof}[Proof of Theorem \upshape\ref{Thm:LDE-long}]
\ref{Thm:LDE-long-item-1}
Assume $\Phi(\lambda)\sim \lambda^\alpha \ell(1/\lambda)$ as $\lambda\to 0+$. 
We use Lemma \ref{Lem:Gt-estimate} with  $c=\varepsilon(t)$ to get
\begin{align}
   U[0,\varepsilon(t)t]\Big({\rm k}+\Pi(t,\infty)\Big)
   &\leq
   \bP\bigg(\frac{G(t)}{t}\leq \varepsilon(t)\bigg)
   \\
   &\leq 
   U[0,\varepsilon(t)t]\Big({\rm k}+\Pi(t-\varepsilon(t)t,\infty)\Big).
   \label{Thm:LDE-long-item-1-proof-2}	
\end{align}
(In the current setting ${\rm k}=0$; however, we preserve the parameter $\rm k$ with a view toward the proof of Theorem \ref{Thm:LDE-short}.)
By Proposition \ref{Prop:RV-long}, we have
\begin{align}
    U[0,\varepsilon(t)t]\Big({\rm k}+\Pi(t,\infty)\Big)
    &\sim
    \frac{(\varepsilon(t)t)^\alpha}{\ell(\varepsilon(t)t)\Gamma(1+\alpha)}
    \cdot
    \frac{t^{-\alpha}\ell(t)}{\Gamma(1-\alpha)}
    \\
    &=
    \frac{\sin(\pi\alpha)}{\pi\alpha}\frac{\ell(t)}{\ell(\varepsilon(t)t)}\varepsilon(t)^\alpha,
    \quad
    \text{as $t\to\infty$.}
    \label{Thm:LDE-long-item-1-proof-3}
\end{align}
Here we used Euler's reflection formula $\Gamma(\alpha)\Gamma(1-\alpha)=\pi/\sin(\pi\alpha)$.
On the other hand, by the uniform convergence theorem for regularly varying functions \cite[Theorem 1.5.2]{BinGolTeu}, we see that
\begin{align}
	{\rm k}+\Pi(t-\varepsilon(t)t,\infty)\sim {\rm k}+\Pi(t,\infty),
	\quad
	\text{as $t\to\infty$}.
	    \label{Thm:LDE-long-item-1-proof-4}
\end{align}
The desired estimate \eqref{Thm:LDE-long-1} follows from \eqref{Thm:LDE-long-item-1-proof-2}, \eqref{Thm:LDE-long-item-1-proof-3} and \eqref{Thm:LDE-long-item-1-proof-4}.

We further assume that $\tau$ is special. We use Lemma \ref{Lem:Gt-estimate} with ${\rm k}=0$ and $c=\varepsilon(t)$  to get
\begin{align}
	&u(t)\bigg({\rm d}+\int_0^{\varepsilon(t)t}\Big({\rm k}+\Pi(x,\infty)\Big)\,dx\bigg)
	\\
	&\leq
	\bP\bigg(\frac{G(t)}{t}\geq 1-\varepsilon(t)\bigg)
	\\
	&\leq
	u(t-\varepsilon(t)t)\bigg({\rm d}+\int_0^{\varepsilon(t)t}\Big({\rm k}+\Pi(x,\infty)\Big)\,dx\bigg).
	\label{Thm:LDE-long-item-1-proof-7}
\end{align}
Proposition \ref{Prop:RV-long} implies
\begin{align}
	u(t)\bigg({\rm d}+\int_0^{\varepsilon(t)t}\Big({\rm k}+\Pi(x,\infty)\Big)\,dx\bigg)
	&\sim
	\frac{t^{-1+\alpha}}{\ell(t)\Gamma(\alpha)}\cdot
    \frac{(\varepsilon(t)t)^{1-\alpha}\ell(\varepsilon(t)t)}{\Gamma(2-\alpha)}
    \\
    &=
    \frac{\sin(\pi \hat{\alpha})}{\pi\hat{\alpha}}\frac{\ell(\varepsilon(t)t)}{\ell(t)}\varepsilon(t)^{\hat{\alpha}},
    \quad
    \text{as $t\to\infty$,}
    \label{Thm:LDE-long-item-1-proof-8}
\end{align}
where $\hat{\alpha}=1-\alpha$.
On the other hand, the uniform convergence theorem  yields that
\begin{align}
    u(t-\varepsilon(t)t)\sim
    u(t),
    \quad \text{as $t\to\infty$.}
    \label{Thm:LDE-long-item-1-proof-9}	
\end{align}
We use \eqref{Thm:LDE-long-item-1-proof-7}, \eqref{Thm:LDE-long-item-1-proof-8} and \eqref{Thm:LDE-long-item-1-proof-9} to obtain \eqref{Thm:LDE-long-2}, as desired.

\ref{Thm:LDE-long-item-2}
Assume $\Phi(\lambda)\asymp \lambda^\alpha \ell(1/\lambda)$ as $\lambda\to 0+$.
We use Proposition \ref{Prop:mild-RV-long} to get
\begin{align}
   U[0,\varepsilon(t)t]\Big({\rm k}+\Pi(t,\infty)\Big)
   \asymp
   \frac{(\varepsilon(t)t)^{\alpha}}{\ell(\varepsilon(t)t)}\cdot t^{-\alpha}\ell(t)=\frac{\ell(t)}{\ell(\varepsilon(t)t)}\varepsilon(t)^\alpha,
   \quad
   \text{as $t\to\infty$.}
   \label{Thm:LDE-long-item-2-proof-1}
\end{align}
In addition, the uniform convergence theorem implies
\begin{align}
  	\frac{{\rm k}+\Pi(t-\varepsilon(t)t,\infty)}{{\rm k}+\Pi(t,\infty)}
  	\asymp
  	\frac{\ell(t-\varepsilon(t)t)}{\ell(t)}\to 1,
  	\quad
  	\text{as $t\to\infty$.}
  	\label{Thm:LDE-long-item-2-proof-2}
\end{align}
Combining \eqref{Thm:LDE-long-item-1-proof-2} with \eqref{Thm:LDE-long-item-2-proof-1} and \eqref{Thm:LDE-long-item-2-proof-2}, we obtain \eqref{Thm:LDE-long-3}. 

We further assume that $\tau$ is special. Then we have \eqref{Thm:LDE-long-item-1-proof-7}. We use Proposition \ref{Prop:mild-RV-long} to get
\begin{align}
	u(t)\bigg({\rm d}+\int_0^{\varepsilon(t)t}\Big({\rm k}+\Pi(x,\infty)\Big)\,dx\bigg)
	&\asymp
	\frac{t^{-1+\alpha}}{\ell(t)}
	\cdot
	(\varepsilon(t)t)^{1-\alpha}\ell(\varepsilon(t)t)
	\\
	&=\frac{\ell(\varepsilon(t)t)}{\ell(t)}\varepsilon(t)^{\hat{\alpha}},
	\quad
	\text{as $t\to\infty$.}
	\label{Thm:LDE-long-item-2-proof-3}
\end{align}
In addition, the uniform convergence theorem implies
\begin{align}
    \frac{u(t-\varepsilon(t)t)}{u(t)}
    \asymp
    \frac{\ell(t-\varepsilon(t)t)}{\ell(t)}
    \to 1,
    \quad
    \text{as $t\to\infty$.}
    \label{Thm:LDE-long-item-2-proof-4}	
\end{align}
Combining \eqref{Thm:LDE-long-item-1-proof-7} with \eqref{Thm:LDE-long-item-2-proof-3} and \eqref{Thm:LDE-long-item-2-proof-4}, we obtain \eqref{Thm:LDE-long-4}. 
\end{proof}

We can also prove Theorem \ref{Thm:LDE-short} by replacing $\lambda\to 0+$ and $t\to\infty$ with $\lambda\to\infty$ and $t\to0+$ respectively in the proof of Theorem \ref{Thm:LDE-long}.

\section{Applications to the Last Visit Times of Markov Processes}\label{Sec:Dual}

In this section we summarize the abstract theory of dual processes, local times and inverse local times. We will see that an inverse local time is a subordinator, and our results are applicable to the last visit time of a suitable class of Markov processes. We refer the reader to Blumenthal and Getoor \cite{BluGet64, BluGet, BluGet70}, Getoor and Sharpe \cite{GetSha}, Chung and Walsh \cite{ChuWal}, Marcus and Rosen \cite{MarRos} for more details. 

Let $(S,d)$ be a locally compact separable metric space. 
We adjoin the element $\partial$ to $S$ as the point at infinity if $S$ is non-compact, and as an isolated point if $S$ is compact. In any case $\partial$ is regarded as a cemetery point.
Denote by $\cB(S)$ the Borel $\sigma$-fields of $S$.
We write $\cB(S)_+$ for the set of non-negative measurable functions $f:S\to[0,\infty]$. For $f\in \cB(S)_+$, we adopt the convention that $f(\partial)=0$. 
Let us impose the following assumption:
\begin{enumerate}[label=\upshape(A\arabic*)]
\item\label{Asmp-Hunt} $X=(\Omega, \cM, \cM_t, X_t,\theta_t,\bP_x)$ and $\hat{X}=(\hat{\Omega}, \hat{\cM}, \hat{\cM}_t,\hat{X}_t,\hat{\theta_t},\hat{\bP}_x)$ are Hunt processes  with common state space $(S,\cB(S))$ in the sense of \cite[Section I.9]{BluGet}.
\end{enumerate}
Denote by $(P_t)_{t\geq0}$ and $(\hat{P}_t)_{t\geq0}$ the families of the transition kernels or operators of $X$ and $\hat{X}$, and by $(R_\lambda)_{\lambda>0}$ and $(\hat{R}_\lambda)_{\lambda>0}$ the resolvents of $(P_t)_{t\geq0}$ and $(\hat{P}_t)_{t\geq0}$, that is, the families of resolvent (or potential) kernels or operators of $X$ and $\hat{X}$, respectively. In other words,
\begin{align}
   P_t(x,B)&=\bP_x(X_t\in B),
   &
   P_t f(x)&=\bE_x f(X_t),
\\
R_\lambda(x,B)&=\int_0^\infty e^{-\lambda t}P_t(x,B)\,dt,
	&
	R_\lambda f(x)&=\int_0^\infty e^{-\lambda t} P_tf(x)\,dt,
\end{align}
for $t\geq0$, $\lambda>0$,  $x\in S$, $B\in \cB(S)$, $f\in \cB(S)_+$. Similarly,
\begin{align}
   \hat{P}_t(x,B)&=\hat \bP_x(\hat{X}_t\in B),
   &
   \hat{P}_t f(x)&=\bE_x f(\hat{X}_t),
\\
\hat{R}_\lambda(x,B)&=\int_0^\infty e^{-\lambda t}\hat{P}_t(x,B)\,dt,
	&
	\hat{R}_\lambda f(x)&=\int_0^\infty e^{-\lambda t} \hat{P}_tf(x)\,dt.
\end{align}
We also assume the following condition:
\begin{enumerate}[resume, label=\upshape(A\arabic*)]
\item\label{Asmp-Dual} Let $m$ be a Radon measure on $(S,\cB(S))$. The resolvents $(R_\lambda)_{\lambda>0}$ and $(\hat{R}_\lambda)_{\lambda>0}$ are in duality relative to $m$. In other words, the measures $R_\lambda(x,\cdot )$ and $\hat{R}_\lambda(x,\cdot)$ are absolutely continuous with respect to $m$ for any $\lambda>0$ and $x\in S$, and $ \int_S f(x)R_\lambda g(x)\,m(dx)
     =
     \int_S \hat{R}_\lambda f(x) g(x)\,m(dx)$ for any $\lambda>0$ and $f, g\in \cB(S)_+$.

\end{enumerate}
For $\lambda>0$,
a function $f\in \cB(S)_+$ is said to be $\lambda$-excessive (respectively, $\lambda$-coexcessive) if $\mu R_{\lambda+\mu}f\leq f$ for all $\mu>0$ and $\mu R_{\lambda+\mu}\uparrow f$ as $\mu\to\infty$ (respectively, $\mu \hat{R}_{\lambda+\mu}f\leq f$ for all $\mu>0$ and $\mu \hat{R}_{\lambda+\mu}\uparrow f$ as $\mu\to\infty$).   
By \cite[Theorem 13.2]{ChuWal} (see also \cite[Theorem (VI.1.4)]{BluGet} and \cite[Section 2]{GetSha}),
for each $\lambda>0$ there exists a measurable function $r_\lambda:S\times S\to [0,\infty]$ such that
\begin{enumerate}[label=(\roman*)]
	\item the function $x\mapsto r_\lambda(x,y)$ is $\lambda$-excessive for any $\lambda>0$ and $y\in S$; 
	\item the function $y\mapsto r_\lambda(x,y)$ is $\lambda$-coexcessive for any $\lambda>0$ and $x\in S$; 
	\item $R_\lambda(x,dy)=r_\lambda (x,y)\,m(dy)$ and $\hat{R}_\lambda(y,dx)=m(dx)\,r_\lambda(x,y)$ for any $\lambda>0$ and  $x,y\in S$.
\end{enumerate}
Such a function $r_\lambda$ is uniquely determined, which is guaranteed by \cite[Proposition (VI.1.3)]{BluGet}.
We call $r_\lambda(x,y)$ the $\lambda$-resolvent (or $\lambda$-potential) density. By the resolvent equation 
\begin{align}
R_{\lambda}-R_{\mu}=(\mu-\lambda)R_{\lambda}R_{\mu},
\quad\lambda,\mu>0,
\end{align}
 and the fact $y\mapsto r_\lambda(x,y)$ is $\lambda$-coexcessive, we have
\begin{align}\label{eq:resolvent-eq}
	r_{\lambda}(x,y)-r_{\mu}(x,y)
	=
	(\mu-\lambda)\int_S r_{\lambda} (x,z)r_{\mu}(z,y)\,m(dz),
\end{align}
for any $\lambda,\mu>0$ and $x,y\in S$, as shown in \cite[Theorem 13.2]{ChuWal}.
By \cite[Remark (VI.1.13)]{BluGet}, the measure $m$ must be a reference measure of $X$ and $\hat{X}$, that is, for any $B\in \cB(S)$, the condition $m(B)=0$  implies that $\int_0^\infty P_t(x,B)\,dt=\int_0^\infty \hat{P}_t(x,B)\,dt=0$ for any $x\in S$. In particular, the topological support of $m$ must coincide with $S$.

For $B\in \cB(S)$, we define $T_B=\inf\{t>0\::\: X_t\in B\}$ and $\hat{T}_B=\inf\{t>0\::\:\hat{X}_t\in B\}$. A point $a\in S$ is said to be regular for $B$ provided that $\bP_{a}(T_B=0)=1$. The following proposition states a sufficient condition for $a$ to be regular for itself. Similar propositions can be found in \cite[Proposition (VI.4.11)]{BluGet} and \cite[Proposition (7.3)]{BluGet70} under slightly different settings.

\begin{Prop}\label{Prop:regular}
Suppose \ref{Asmp-Hunt} and \ref{Asmp-Dual} are satisfied.
Let $\lambda>0$ and $a\in S$. Assume that $r_\lambda(a,a)<\infty$ and  the function $x\mapsto r_\lambda(x,a)$ is continuous at $x=a$. Then $r_\lambda(a,a)>0$ and
\begin{align}
     \bE_x\exp(-\lambda T_{\{a\}})
     =
     \frac{r_\lambda(x,a)}{r_\lambda(a,a)},
     \quad
     x\in S.
     \label{Prop:regular-1}
\end{align}
In particular, $a$ is regular for $\{a\}$. 
\end{Prop}

\begin{proof}

By the switching formula \cite[Theorem (VI.1.16)]{BluGet},
\begin{align}
    \bE_x [\exp(-\lambda T_B)r_\lambda (X_{T_B},a)\:;\: T_B<\infty]
    =
    \hat{\bE}_{a} [\exp(-\lambda \hat{T}_B)
    r_\lambda(x, \hat{X}_{\hat{T}_B})\:;\:\hat{T}_B<\infty], 
    \label{Prop:regular-proof-1}
\end{align}
for any $x\in S$ and $B\in \cB(S)$.
Let $\{B_n\}_{n=1}^\infty\subset S$ be a sequence of precompact open sets with 
$B_n\supset B_{n+1}$ ($n=1,2,\dots$) and $\bigcap_{n=1}^{\infty} \bar{B_n}=\bigcap_{n=1}^\infty B_n=\{a\}$. 
Then $T:=\sup_{n}T_{B_n}\,(\leq T_{\{a\}})$ coincides with $T_{\{a\}}$, $\bP_x$-a.s. (Indeed, we use the quasi-left continuity of $(X_t)$ to see that 
\begin{align}
X_{T}=\lim_{n\to\infty} X_{T_{B_n}}\in \bigcap_{n=1}^\infty \bar{B_n}=\{a\},
\quad\text{$\bP_x$-a.s.\ on $\{T<\infty\}$},
\end{align} 
 and hence $T=T_{\{a\}}$, $\bP_x$-a.s.\ on $\{T<\infty\}$. Here the precompactness of $B_n$ ensured that $X_T\neq \partial$, $\bP_x$-a.s.\ on $\{T<\infty\}$. In addition, $T_{\{a\}}=\infty$ on $\{T=\infty\}$.) 
By the right-continuity of $(\hat{X}_t)$, we have $\hat{T}_{B_n}=0$, $\hat{\bP}_a$-a.s.
Hence we substitute $B=B_n$ to get
\begin{align}\label{Prop:regular-proof-2}
    \bE_x[\exp(-\lambda T_{B_n})r_\lambda(X_{T_{B_n}},a)\:;\:T_{B_n}<\infty]
       =
       r_\lambda(x,a). 
\end{align}
 Since $T_{B_n}\uparrow T_{\{a\}}$, $\bP_x$-a.s., $r_\lambda(a,a)<\infty$ and the function $x\mapsto r_\lambda(x,a)$ is continuous at $x=a$, we use   the dominated convergence theorem and \eqref{Prop:regular-proof-2} to get
\begin{align}
	\bE_x\exp(-\lambda T_{\{a\}})r_\lambda(a,a)
	&=\lim_{n\to\infty}\bE_x[\exp(-\lambda T_{B_n})r_\lambda(X_{T_{B_n}},a)\:;\:T_{B_n}<\infty]
	\\
	&=r_\lambda(x,a).
\end{align}
In particular, $r_\lambda(a,a)\geq r_\lambda(x,a)$.
Since the function $x\mapsto r_\lambda(x,a)=\hat{R}_\lambda(a,dx)/m(dx)$ is not identically zero, $r_\lambda(a, a)$ must be positive, and it is finite by the assumption.
Therefore we obtain \eqref{Prop:regular-1}. By substituting $x=a$ into \eqref{Prop:regular-1}, we have $\bE_{a}\exp(-\lambda T_{\{a\}})=1$. Hence $a$ is regular for $\{a\}$.
\end{proof}

Let $a\in S$. From now on, we suppose
\begin{enumerate}[label=\upshape(A\arabic*)]
\setcounter{enumi}{2}
\item\label{Asmp-u1} $r_1(a,a)<\infty$ and  the function $x\mapsto r_1(x,a)$ is continuous at $x=a$.	
\end{enumerate}
By Proposition \ref{Prop:regular} and \cite[Theorem (V.3.13)]{BluGet}, there exists a local time $L$ at $a$. It is a perfect positive continuous additive functional with $\mathop{\rm Supp} L=\{a\}$,
which is uniquely determined up to a multiplicative constant.
Its
$1$-potential is given by $\bE_x\int_0^\infty e^{-t}\,dL(t) = c\bE_x \exp(-T_{\{a\}})$, $x\in S$, for some constant $c>0$. 
\begin{Prop}\label{Prop:normalization-local-time}
Suppose \ref{Asmp-Hunt}, \ref{Asmp-Dual}, \ref{Asmp-u1} are satisfied.
Then the local time $L$ at $a$ can be normalized so that
\begin{align}
	\bE_x\int_0^\infty e^{-\lambda t}\,dL(t)
	=
	r_\lambda(x,a),
	\quad
	\lambda>0,\,x\in S.
\end{align}
\end{Prop}

\begin{proof}
Let us choose the normalization $c=r_1(a,a)$ for the local time $L$. Then Proposition \ref{Prop:regular} implies that  $\bE_x\int_0^\infty e^{-t}\,dL(t)= r_1(x,a)$, $x\in S$. We use \cite[Proposition (IV.2.4)]{BluGet} and \eqref{eq:resolvent-eq} to get
\begin{align}
    \bE_x\int_0^\infty e^{-\lambda t}\,dL(t) 
    &= r_1(x,a)
    +(1-\lambda)
    \int_S r_\lambda(x,y) r_1(y,a)\,m(dy)
    \\
    &=r_\lambda(x,a),
    \qquad
    \lambda>0,
    \,x\in S,
\end{align}
as desired.	
\end{proof}

Throughout the rest of this section we normalize $L$ as in Proposition \ref{Prop:normalization-local-time}.
Let $\tau$ be the inverse local time at $a$, that is,
\begin{align}
\tau(s)=\inf\{t>0\::\: L(t)>s\}, \quad s\geq0.
\end{align}
Here it is understood that $\tau(s)=\infty$ when $s\geq L(\infty):=\lim_{t\to\infty }L(t)$.

\begin{Prop}\label{Prop:local-time-resolvent}
Suppose \ref{Asmp-Hunt}, \ref{Asmp-Dual}, \ref{Asmp-u1} are satisfied.
Then the inverse local time $(\tau, \bP_a)$ is a subordinator with Laplace exponent
\begin{align}
   \Phi(\lambda)
   =
   \bigg(\bE_{a}\int_0^\infty e^{-\lambda t}\,dL(t)\bigg)^{-1}
   =
   \frac{1}{r_\lambda(a,a)},
   \quad
   \lambda>0.	
\end{align}
 In addition,
\begin{align}
   \bar{\{\tau(s)\::\:0\leq s< L(\infty)\}}
   =
   \bar{\{t\geq0\::\: X_t=a\}},
   \quad
   \text{$\bP_{a}$-a.s.}	
\end{align}	
\end{Prop}

\begin{proof}
This is a consequence of 	\cite[Theorems (V.3.21) and (V.3.8)]{BluGet}. See also \cite[(2.4)]{BluGet64}.
\end{proof}

Set
\begin{align}
   G(t)=\sup\{\tau(s)\::\: \text{$s\geq0$, $\tau(s)\leq t$}\}=
   \sup\{0\leq s\leq t\::\:X_s=a\},
   \quad
   t\geq0,
   \,\text{$\bP_{a}$-a.s.}	
\end{align}
We use Theorems \ref{Thm:LDE-long} and \ref{Thm:LDE-short} to obtain the following result:

\begin{Thm}[Long-range large deviation estimates]\label{Thm:LV-LDE-long}
Suppose \ref{Asmp-Hunt}, \ref{Asmp-Dual}, \ref{Asmp-u1} are satisfied. Let $\alpha\in(0,1)$ and $\hat{\alpha}=1-\alpha$. Let $\ell:(0,\infty)\to(0,\infty)$ be slowly varying at $\infty$, and suppose $\varepsilon:(0,\infty)\to(0,1)$ satisfies $\varepsilon(t)\to 0$ and $\varepsilon(t)t\to\infty$ as $t\to\infty$.
 Then the following statements hold.
 \begin{enumerate}[label=\upshape(\arabic*)]
\item Assume $r_\lambda(a,a)\sim \lambda^{-\alpha}/ \ell(1/\lambda)$ as $\lambda\to0+$. Then
\begin{align}\label{Thm:LV-LDE-long-1}
  	 \mathbb{P}_{a}\bigg(\frac{G(t)}{t}\leq \varepsilon(t)\bigg)
     \sim
     \frac{\sin(\pi\alpha)}{\pi\alpha}
     \frac{\ell(t)}{\ell(\varepsilon(t)t)}
    \varepsilon(t)^\alpha,
\end{align}
as $t\to\infty$. Suppose additionally that $\tau$ is special. Then
\begin{align}\label{Thm:LV-LDE-long-2}
	 \mathbb{P}_{a}\bigg(\frac{G(t)}{t}\geq 1- \varepsilon(t)\bigg)
     \sim
     \frac{\sin(\pi\hat{\alpha})}{\pi\hat{\alpha}}
     \frac{\ell(\varepsilon(t)t)}{\ell(t)}
     \varepsilon(t)^{\hat{\alpha}},
\end{align}
as $t\to\infty$.

\item Assume $r_\lambda(a, a)\asymp \lambda^{-\alpha}/\ell(1/\lambda)$ as $\lambda\to0+$. Then
\begin{align}\label{Thm:LV-LDE-long-3}
  	 \mathbb{P}_{a}\bigg(\frac{G(t)}{t}\leq \varepsilon(t)\bigg)
     \asymp
     \frac{\ell(t)}{\ell(\varepsilon(t)t)}
     \varepsilon(t)^\alpha,
\end{align}
as $t\to\infty$.
Suppose additionally that $\tau$ is special. Then
\begin{align}\label{Thm:LV-LDE-long-4}
	 \mathbb{P}_{a}\bigg(\frac{G(t)}{t}\geq 1- \varepsilon(t)\bigg)
     \asymp
     \frac{\ell(\varepsilon(t)t)}{\ell(t)}
     \varepsilon(t)^{\hat{\alpha}},
\end{align}
as $t\to\infty$.

\end{enumerate}
\end{Thm}

\begin{Thm}[Short-range large deviation estimates]\label{Thm:LV-LDE-short}
Suppose \ref{Asmp-Hunt}, \ref{Asmp-Dual}, \ref{Asmp-u1} are satisfied. 
Let $\alpha\in(0,1)$  and $\hat{\alpha}=1-\alpha$. Let $\ell:(0,\infty)\to(0,\infty)$ be slowly varying at $0+$, and suppose $\varepsilon:(0,\infty)\to(0,1)$ satisfies $\varepsilon(t)\to 0$ as $t\to0+$.
Then the following statements hold.
\begin{enumerate}[label=\upshape(\arabic*)]
\item Assume  $r_\lambda(a,a)\sim \lambda^{-\alpha}/ \ell(1/\lambda)$ as $\lambda\to\infty$, then \eqref{Thm:LV-LDE-long-1} holds as $t\to0+$. Suppose additionally that $\tau$ is special. Then
\eqref{Thm:LV-LDE-long-2} holds as $t\to0+$.

\item Assume  $r_\lambda(a,a)\asymp \lambda^{-\alpha}/ \ell(1/\lambda)$ as $\lambda\to\infty$, then \eqref{Thm:LV-LDE-long-3} holds as $t\to0+$. Suppose additionally that $\tau$ is special. Then
\eqref{Thm:LV-LDE-long-4} holds as $t\to0+$.

\end{enumerate}
\end{Thm}

To conclude this section, we provide sufficient conditions for $(\tau, \bP_{a})$ to be a special or complete subordinator.
By \cite[Proposition 13.6]{ChuWal}, the semigroups $(P_t)$ and $(\hat{P}_t)$ are in weak duality relative to $m$, that is, $\int_S f(x) P_t g(x)\,m(dx)=\int_S \hat{P}_t f(x)g(x)\,m(dx)$, for any $t>0$ and $f,g\in \cB(S)_+$.
We further assume that
\begin{enumerate}[label=\upshape(A\arabic*)]
	\setcounter{enumi}{3}
	\item\label{Asmp-Dual-P} $P_t(x,\cdot)$ and $\hat{P}_t(y,\cdot)$ are absolutely continuous with respect to $m$ for any $t>0$ and $x, y\in S$.
\end{enumerate}
Then, by \cite[Theorem 1]{Wit}, there  uniquely exists a measurable function $(0,\infty)\times S\times S\ni(t,x,y) \mapsto p_t(x,y)\in[0,\infty]$ such that 
\begin{enumerate}[label=(\roman*)]
	\item $P_t (x,dy)= p_t(x,y)\,m(dy)$ for any $t>0$ and $x,y\in S$;
	\item  $\hat{P}_t (y, dx)= p_t(x,y)\,m(dx)$ for any $t>0$ and $x,y\in S$;
	\item (Chapman--Kolmogorov equation) $p_{s+t}(x,y)=\int_S p_s(x,z)p_t(z,y)\,m(dz)$ for any $s,t>0$ and $x,y\in S$.
\end{enumerate}

\begin{Prop}\label{Prop:transition-potential}
Suppose \ref{Asmp-Hunt}, \ref{Asmp-Dual}, \ref{Asmp-u1}, \ref{Asmp-Dual-P} are satisfied. Then
\begin{align}\label{Prop:transition-potential-1}
\frac{1}{\Phi(\lambda)}=r_\lambda(a,a)=\int_0^\infty e^{-\lambda t} p_t(a,a)\,dt,
\quad
\lambda>0.
\end{align}
In other words, the function $t\mapsto p_t(a,a)$ is a density of the potential measure of the subordinator $(\tau,\bP_a)$.
In addition, if $t\mapsto p_t(a,a)$ is non-increasing (respectively, completely monotone), then $(\tau, \bP_{a})$ is special (respectively, complete). 	
\end{Prop}

\begin{proof}
As shown in \cite[Section 3]{GetSha},
\begin{align}
    r_\lambda(x,y)=\int_0^\infty e^{-\lambda t}p_t(x,y)\,dt,
    \quad
    \lambda>0,\, x,y\in S,
\end{align}	
which immediately implies \eqref{Prop:transition-potential-1}. The remainder of the claim follows from Proposition \ref{Prop:special-complete}.
\end{proof}

Let us impose
\begin{enumerate}[label=\upshape(A\arabic*)]
	\setcounter{enumi}{4}
		\item\label{Asmp-sym} $X$ is $m$-symmetric, that is, $X=\hat{X}$. 
\end{enumerate}
Then $p_t(x,y)=p_t(y,x)$ for any $t>0$ and $x,y\in S$. We denote by $\|\cdot\|$ the norm on $L^2(m)$. By Jensen's inequality, $|P_t f(x)|^2\leq P_t(|f|^2)(x)$ ($t>0$, $f\in L^2(m)$, $x\in S$) and hence
\begin{align}
   \|P_t f\|^2
   \leq
   \int_S P_t(|f|^2)(x)\,m(dx)
   =
   \int_S |f(x)|^2 P_t1(x)\,m(dx)
   \leq \|f\|^2,
   \quad
   t>0,\,f\in L^2(m). 	
\end{align}
Here we used $P_t 1(x)\leq 1$ for $t\geq0$ and $x\in S$.
Hence $P_t$ is a contraction on $L^2(m)$. We can also see that $P_t$ is self-adjoint on $L^2(m)$. Therefore, for any $0<s<t$,
\begin{align}
   p_t(a,a)=\|p_{t/2}(a,\cdot)\|^2=\| P_{(t-s)/2}p_{s/2}(a,\cdot )\|^2\leq \|p_{s/2}(a,\cdot)\|^2=p_s(a,a).
\end{align}
In other words, the map $t\mapsto p_t(a,a)$ is non-increasing. By Proposition \ref{Prop:transition-potential}, we conclude that $(\tau, \bP_a)$ is a special subordinator.
Finally we also consider the following condition.

\begin{enumerate}[label=\upshape(A\arabic*)]
	\setcounter{enumi}{5}
	\item\label{Asmp-sym-strong} $(P_t)$ is a  strongly continuous semigroup on $L^2(m)$.		
\end{enumerate}

\begin{Prop}\label{Prop:symmetric-complete}
Assume \ref{Asmp-Hunt}, \ref{Asmp-Dual}, \ref{Asmp-u1}, \ref{Asmp-Dual-P}, \ref{Asmp-sym}, \ref{Asmp-sym-strong} are satisfied. Then  the map $t\mapsto p_t(a,a)$ is completely monotone, and hence $(\tau, \bP_a)$ is a complete subordinator.
\end{Prop}

\begin{proof}
Let $-A$ be the $L^2(m)$-generator of $(P_t)$, that is, $P_t=e^{-tA}$, $t\geq0$. The spectrum of $P_t$ is contained in $[0,1]$, and hence that of $A$ is contained in $[0,\infty)$.
We denote by $E$ the spectral decomposition of $A$. Then $P_t=\int_{[0,\infty)} e^{-t\lambda}\,dE(\lambda)$, $t\geq0$ by the spectral theorem.  See for example \cite[Theorem 13.38 and Note below]{Rud-FA}. Hence
\begin{align}
   \langle P_t f, g\rangle =\int_{[0,\infty)} e^{-t\lambda}\,d\langle E(\lambda) f,g \rangle,	
   \quad
   t\geq 0,\,\,f,g\in L^2(m),
\end{align}
where $\langle\cdot, \cdot\rangle$ denotes the inner product on $L^2(m)$. Fix $\varepsilon>0$ arbitrarily. Then
\begin{align}
  p_{t}(a,a)
  &=
  \big\langle p_{t-\varepsilon}(a,\cdot), p_{\varepsilon}(a,\cdot)\big\rangle
  =
  \big\langle P_{t-2\varepsilon}p_{\varepsilon}(a,\cdot), p_{\varepsilon}(a,\cdot)\big\rangle
  \\
  &=\int_{[0,\infty)} e^{-(t-2\varepsilon)\lambda}\,d\big\langle E(\lambda) p_\varepsilon(a,\cdot), p_\varepsilon(a,\cdot)\big\rangle,
   \quad
   t\geq 2\varepsilon.
   \label{p_t-CM}
\end{align}
Note that the map $\lambda\mapsto \langle E(\lambda) f, f\rangle$ is non-negative and non-decreasing with $\langle E(\lambda) f, f\rangle\to \|f\|^2 <\infty$ as $\lambda\to\infty$ for each $f\in L^2(m)$. By \eqref{p_t-CM}, we see that
$(-\partial /\partial t)^n p_t(a,a)\geq0$ for any  $n\in \bN$ and $t>2\varepsilon$.
Since $\varepsilon>0$ was arbitrary, the map $t\mapsto p_t(a,a)$ is completely monotone.
The remainder of the claim follows from Proposition \ref{Prop:transition-potential}. 
\end{proof}

\section{Examples: One-dimensional Generalized Diffusion Processes}\label{Sec:Diffusion}

In this section we recall basic notions of generalized diffusion processes (which are also called gap diffusions or quasi diffusions), following
Watanabe \cite{Wat74, Wat95}, Kasahara \cite{Kas}, Knight \cite{Kni}, Kotani and Watanabe \cite{KotWat}, Bertoin \cite[Chapter 9]{Ber99}, Schilling, Song, and Vondra\v{c}ek  \cite[Chapter 15]{SchSonVon}. 
We also refer the reader to Dym and McKean \cite{DymMcK} for the details of Krein's string theory.
Then we apply our abstract results in Section \ref{Sec:Dual} to generalized diffusion processes.

A non-decreasing and right-continuous function $m:[0,\infty]\to[0,\infty]$ with $m(\infty)=\infty$ is called a \emph{string}. The family of all strings will be denoted by $\cM$. For $m\in \cM\setminus\{\infty\}$, set $m(0-)=0$ and $\ell=\sup\{x\geq0\::\: m(x)<\infty\}\in(0,\infty]$.  Then $m(\ell)=\infty$ by the definition. The string $m\in \cM\setminus\{\infty\}$ can be identified as the Lebesgue--Stieljes measure on $[0,\ell)$ associated with $m$, which will be denoted by $m(dx)$.

Fix $m\in \cM\setminus\{\infty\}$ from now on. Denote by $S$ the support of the measure $m(dx)$ in $[0,\ell)$. Let $((B(t))_{t\geq0}, (\bP_x)_{x\in [0,\ell)})$ be a one-dimensional  Brownian motion reflected at $0$ and killed at the height $\ell$ with the cemetery point $\partial$. We denote by $(L^B(t,x)\::\:t\geq0, \, x\in[0,\ell))$
the jointly continuous version of local times of $(B(t))_{t\geq0}$: 
\begin{align}
    \frac{1}{2}\int_0^t f(B(s))\,ds
    =
    \int_0^\ell f(x)L^B(t,x)\,dx,
    \quad
    t\geq0,
    \label{eq:Brownian-LT}
\end{align}
for any non-negative Borel function $f:\bR\to [0,\infty)$. Here it is understood that $f(\partial)=0$. 
We denote by $(A(t))_{t\geq0}$  the continuous additive functional associated with $m(dx)$:
\begin{align}
    A(t)=\int_S L^B(t,x)\,m(dx),
    \quad
    t\geq0.
\end{align}
Set $A^{-1}(t)=\inf\{s>0\::\: A(s)>t\}$. We define $(X_t)_{t\geq0}$ as the time changed process
\begin{align}
	X_t=
	\begin{cases}
		B(A^{-1}(t)),
		& t< A(\infty),
		\\
		\partial,
		& t\geq A(\infty).
	\end{cases}
\end{align}
We call $((X_t)_{t\geq0}, (\bP_x)_{x\in S})$ the \emph{generalized diffusion process} associated with $m$.
If $m(0)=0$ and $S=[0,\infty)$, then $((X_t)_{t\geq0}, (\bP_x)_{x\geq0})$ is a $[0,\infty)$-valued diffusion process reflected at $0$ with natural scale and speed measure $m(dx)$. 
It is known that the generalized diffusion process satisfies the conditions \ref{Asmp-Hunt}, \ref{Asmp-Dual}, \ref{Asmp-u1}, \ref{Asmp-Dual-P}, \ref{Asmp-sym}, \ref{Asmp-sym-strong} of Section \ref{Sec:Dual} for any $a\in S$. 
See for example \cite[Section 3]{Wat74} and \cite[Theorem 1.1]{Kni}.
In particular, $p_t(x,y)=p_t(y, x)$ and $r_\lambda(x,y)=r_\lambda(y,x)$ for all $t,\lambda>0$ and $x,y\in S$. We can apply Proposition \ref{Prop:symmetric-complete} to this setting. Nevertheless, a more specific result is known for the inverse local time at $0$ when $0\in S$, as we shall see below.

For $m\in \cM\setminus\{0\}$ with $0\in S$, the function $h(\lambda)=r_{\lambda}(0,0)$ ($\lambda>0$) is called the \emph{characteristic function} of $m$. 
By \cite[Theorem 2]{Kni} and \cite[Theorem 7.3]{SchSonVon}, the function 
 $1/h(\lambda)$ is a non-zero complete Bernstein function with $1/h(\infty)=\infty$, or equivalently,
$h(\lambda) \not\equiv 0$ is a Stieltjes function with $h(\infty)=0$ in the sense of \cite[Chapter 2]{SchSonVon}. In other words,
there uniquely exists a measure $\sigma(dt)$ on $[0,\infty)$ such that $\int_{[0,\infty)}\sigma(dt)/(1+t)\in(0,\infty)$ and
\begin{align}
	h(\lambda)
	=
	\int_{[0,\infty)}\frac{\sigma(dt)}{\lambda +t},
	\quad
	\lambda>0.
	\label{eq:Krein}
\end{align}
The measure $\sigma(dt)$ is called the \emph{spectral measure} of $m$.
Conversely, for any Stieltjes function $h(\lambda)\not\equiv 0$ with $h(\infty)=0$, there uniquely exists a string $m\in\cM\setminus\{\infty\}$ such that $0\in S$ and $r_\lambda(0,0)=h(\lambda)$ holds for any $\lambda>0$. The correspondence between $m$ and $h$ is called \emph{Krein's correspondence}.

Throughout the rest of this section, we assume $0\in S$.
Let $(\tau, \bP_0)$ be the inverse local time at $0$, which is a complete (and hence special) subordinator with Laplace exponent $\Phi(\lambda)=1/h(\lambda)$. 
Set $u(t)=\inf\{y>0\::\:x m(x)>t\}$ ($t>0$). By \cite[Theorem 2.3]{KotWat}, we have $h(\lambda)\asymp u(1/\lambda)$, for all $\lambda>0$. In particular, we obtain the following equivalence of asymptotic behaviors of $m(x)$ and $h(\lambda)$.

\begin{Prop}\label{Prop:BRV-m-h}
Let $0<\alpha<1$.
Let $f:(0,\infty)\to(0,\infty)$ be a regularly varying function at $\infty$ (respectively, at $0+$) with index $1/\alpha$. 
Take a measurable function $g:(0,\infty)\to (0,\infty)$ so that $f(g(x))\sim g(f(x))\sim x$ as $x\to \infty$ (respectively, as $x\to0+$).
Then the following conditions are equivalent:
\begin{enumerate}[label=\upshape(\arabic*)]
	\item $m(x)\asymp f(x)/x$, as $x\to\infty$ (respectively, as $x\to 0+$).
	\item $h(\lambda) \asymp g(1/\lambda)$, as $\lambda\to 0+$ (respectively, as $\lambda\to\infty$).
\end{enumerate}
\end{Prop}
Note that $g(1/\lambda)$ is regularly varying at $0+$ (respectively, at $\infty$) with index $-\alpha$ by \cite[Theorem 1.5.12]{BinGolTeu}.
The following theorem is a sharper result than the previous proposition:

\begin{Thm}[{\cite[Theorem 2]{Kas}, \cite[Corollary of Theorem 2.1]{KotWat}}]\label{Thm:RV-m-h}
Let $0<\alpha<1$. Let $f$ and $g$ be as in Proposition {\upshape\ref{Prop:BRV-m-h}}. Set
\begin{align}
	D_\alpha
	=
	\frac{\Gamma(1+\alpha)}{\Gamma(1-\alpha)}(\alpha(1-\alpha))^{-\alpha}.
\end{align}
Then the following conditions are equivalent:
\begin{enumerate}[label=\upshape(\arabic*)]
	\item $m(x)\sim f(x)/x$, as $x\to\infty$ (respectively, as $x\to 0+$).
	\item $h(\lambda) \sim D_\alpha g(1/\lambda)$, as $\lambda\to 0+$ (respectively, as $\lambda\to\infty$).
\end{enumerate}
\end{Thm}

\begin{Ex}[Bessel diffusion process]
For $c>0$ and $0<\alpha<1$, we define the string $m$ by
\begin{align}
    m(x)= cx^{(1/\alpha)-1},
    \quad x\geq0.	
\end{align}
Then the corresponding diffusion process $X$ coincides with a $(2-2\alpha)$-dimensional Bessel diffusion process in natural scale up to a multiplicative constant, reflected at $0$.   
In this case, we have $h(\lambda)=c^{-\alpha} D_\alpha \lambda^{-\alpha}$. In other words, $(\tau, \bP_0)$ is a stable subordinator. See Rogers and Williams \cite[V. (48.6)]{RogWil2} and Donati-Martin, Roynette, Vallois and Yor \cite{DMRVY} for more details of Bessel diffusion processes.
\end{Ex}

Therefore we can apply Theorems \ref{Thm:LV-LDE-long} and \ref{Thm:LV-LDE-short} to generalized diffusion processes with $a=0\in S$ and $r_\lambda(0,0)=h(\lambda)$ under the conditions of Proposition \ref{Prop:BRV-m-h} or Theorem \ref{Thm:RV-m-h}.

\section{Examples: One-dimensional L\'evy Processes}\label{Sec:Levy}
In this section we  apply our abstract results in Section \ref{Sec:Dual} to the last visit times of L\'evy processes. We refer the reader to  Bertoin \cite[Section V]{Ber} for more details about local times of  L\'evy processes.

Let $(X_t)_{t\geq0}$ be a real-valued L\'evy process.  
The characteristic function of $X_t$ can be expressed in the form
\begin{align}
    \bE_0 \exp(i\xi X_t)
    =
    \exp(-\Psi(\xi)t),
    \quad
    \xi\in \bR,\,\, t\geq0,	
\end{align}
where $i=\sqrt{-1}$, and  $\Psi:\bR\to \bC$ is called the characteristic exponent of $X$. Note that $\Real \Psi(\xi)\geq0$, $\Psi(-\xi)=\bar{\Psi(\xi)}$ and $\bE_x\exp(i\xi X_t)=\bE_0 \exp(i\xi(X_t+x))$.

From now on we assume that
\begin{align}
 \bigg(\int_{-\infty}^\infty \frac{d\xi}{|\lambda+\Psi(\xi)|}\leq \bigg) 
    \int_{-\infty}^\infty \frac{d\xi}{\lambda+\Real \Psi(\xi)}
    <\infty,
    \quad
    \lambda>0.	
    \label{Asmp:Levy-integrability}
\end{align}	

\begin{Lem}
Suppose that \eqref{Asmp:Levy-integrability} is satisfied. Then the conditions \ref{Asmp-Hunt}, \ref{Asmp-Dual}, \ref{Asmp-u1}, \ref{Asmp-Dual-P} in Section {\upshape\ref{Sec:Dual}} are satisfied for any $a\in\bR$ where $S=\bR$, $\hat{X}=-X$ and $m$ is the Lebesgue measure on $\bR$.	
\end{Lem}

\begin{proof}
It is well known that the condition \ref{Asmp-Hunt} is satisfied.
Note that
\begin{align}
    |\bE_x\exp(i\xi X_t)|
    =
    \exp(-\Real \Psi(\xi) t)
    \leq \frac{1}{1+\Real \Psi(\xi) t}.
\end{align}
By the Fourier inversion (see for example \cite[Theorem 3.3.14]{Dur}), the process $X$ admits a probability density function $p_t(x,y)$ given by
\begin{align}
   p_t(x,y)=p_t(y-x)
   =
   \frac{1}{2\pi}\int_{-\infty}^\infty \exp\Big(-\Psi(\xi)t-i\xi(y-x)\Big)\,d\xi,
\end{align}
and hence the corresponding $\lambda$-resolvent density is given by
\begin{align}
    r_\lambda(x,y)
    =
    r_\lambda(y-x)
    =
    \int_0^\infty e^{-\lambda t}p_t(y-x)\,dt
    =
    \frac{1}{2\pi}\int_{-\infty}^\infty\frac{e^{-i\xi(y-x)}}{\lambda+\Psi(\xi)}\,d\xi.	
\end{align}
It is easily seen that the maps $(0,\infty)\times \bR^2\ni(t,x,y)\mapsto p_t(x,y)\in[0,\infty)$ and $(0,\infty)\times \bR^2\ni(\lambda, x,y)\mapsto r_\lambda(x,y)\in[0,\infty)$ are jointly continuous. 
The probability density function and the $\lambda$-resolvent density of $\hat{X}=-X$ are given by $p_t(x-y)$ and $r_\lambda(x-y)$, respectively.
Therefore \ref{Asmp-Dual}, \ref{Asmp-u1}, \ref{Asmp-Dual-P} are fulfilled.
\end{proof}

Note that $r_\lambda(a,a)=r_\lambda(0)$ can be represented in the form
\begin{align}
    r_\lambda(0)
    =
    \frac{1}{2\pi}\int_{-\infty}^\infty \frac{d\xi}{\lambda+\Psi(\xi)}
    =
    \frac{1}{\pi}\int_{0}^\infty \Real \bigg(\frac{1}{\lambda+\Psi(\xi)}\bigg)\,d\xi,
    \quad
    \lambda>0.
\end{align}
As shown in Proposition \ref{Prop:local-time-resolvent}, the map $\lambda\mapsto 1/r_\lambda(0)$ is the Laplace exponent of the inverse local time at $a$ under $\bP_a$.
\begin{Lem}
Suppose \eqref{Asmp:Levy-integrability} is satisfied.
If $X$ is symmetric, then the map $(0,\infty)\ni t\mapsto p_t(0)\in (0,\infty)$ is completely monotone, and hence
 the map $(0,\infty)\ni\lambda\mapsto 1/r_\lambda(0)\in(0,\infty)$ is a complete (and hence special) Bernstein function.
\end{Lem}

\begin{proof}
We can prove this lemma by checking the conditions of Proposition \ref{Prop:symmetric-complete}. Nevertheless we give a more direct proof here. 
Since $X$ is symmetric, we have $\Psi(\xi)=\Psi(-\xi)\in[0,\infty)$ and hence
\begin{align}
    p_t(0)=\frac{1}{\pi}\int_0^\infty \exp(-\Psi(\xi)t)\,d\xi.  	
\end{align}
which is
 completely monotone. The rest of the claim follows from Proposition \ref{Prop:transition-potential}.
\end{proof}

The following two lemmas give sufficient conditions for the function $\lambda \mapsto r_\lambda(0)$ to be regularly varying or bounded above and below by a regularly varying function, respectively. Therefore we can apply Theorems \ref{Thm:LV-LDE-long} and \ref{Thm:LV-LDE-short} to these settings.

\begin{Lem}[{\cite[Lemmas 4 and 5]{Ber95-FM}}]\label{Lem:RV-levy-local}
Suppose \eqref{Asmp:Levy-integrability} is satisfied.
Let $\beta\in(1,2]$ and $\theta\in \bR$ with $|\theta|\leq (1-\beta/2)\pi$. Assume that the following conditions are fulfilled.
\begin{enumerate}[label=\upshape(\roman*)]
\item\label{Lem:RV-levy-local-1} $|\Psi(\xi)|$ is regularly varying at $\infty$ (respectively, at $0+$) with index $\beta$.
\item	\label{Lem:RV-levy-local-2} $\Psi(\xi)/|\Psi(\xi)|\to e^{i\theta}$, as $\xi\to \infty$ (respectively, as $\xi\to 0+$).
\end{enumerate}
Then the function $\lambda\mapsto r_\lambda(0)$ is regularly varying at $\infty$ (respectively, at $0+$) with index $-\alpha$, where $\alpha=1-1/\beta\in(0,1/2]$. More precisely,
\begin{align}
     r_\lambda(0)
     \sim
     \frac{\cos(\theta/\beta)}{\beta \sin(\pi/\beta)}\frac{|\Psi|^{-1}(\lambda)}{\lambda},
     \quad
     \text{as $\lambda\to\infty$ (respectively, as $\lambda\to 0+$), }
\end{align}
where $|\Psi|^{-1}(\lambda)=\inf\{\xi>0\::\: |\Psi(\xi)|>\lambda\}$ $(\lambda>0)$ is regularly varying at $\infty$ (respectively, at $0+$) with index $1/\beta$.
\end{Lem}

\begin{Rem}
The conditions \ref{Lem:RV-levy-local-1} and \ref{Lem:RV-levy-local-2} in the preceding lemma hold if and only if $(\xi X_{1/|\Psi(\xi)|})_{\xi>0}$ converges in distribution as $\xi\to\infty$ (respectively, as $\xi\to 0+$) to the $\beta$-stable law whose characteristic exponent $\Psi_{\beta,\theta}(\xi)$ is given by
\begin{align}
	\Psi_{\beta,\theta}(\xi)
	=
	|\xi|^\beta \exp(i\theta \operatorname{sgn}\xi),
	\quad
	\xi\in \bR.
\end{align}

\end{Rem}

\begin{Lem}\label{Lem:BRV-levy-local}
Suppose \eqref{Asmp:Levy-integrability} is satisfied.
Let $\beta\in(1,2]$ and let $f:(0,\infty)\to (0,\infty)$ be a function that is regularly varying at $\infty$ (respectively, at $0$) with index $\beta$. Assume
\begin{align}
    \Real \Psi(\xi)\asymp f(\xi) \quad
    \text{and}
    \quad
    \Ima \Psi(\xi)=O(f(\xi)),
    \quad
    \text{as $\xi\to \infty$ (respectively, as $\xi\to0+$).} 	
\end{align}
Then
\begin{align}
      r_\lambda(0)\asymp \frac{f^{-1}(\lambda)}{\lambda}	,
      \quad\text{as $\lambda\to \infty$ (respectively, as $\lambda\to 0+$),}
\end{align}
where $f^{-1}(\lambda)=\inf\{\xi>0\::\: f(\xi)>\lambda\}$ $(\lambda>0)$ is regularly varying at $\infty$ (respectively, at $0+$) with index $1/\beta$.
\end{Lem}

\begin{proof}
We consider only the regime $\xi,\lambda\to\infty$.
Note that
\begin{align}
      \Real \bigg(\frac{1}{\lambda+\Psi(\xi)}\bigg)
      =
      \frac{\lambda+\Real \Psi(\xi)}{(\lambda+\Real \Psi(\xi))^2
     +
     (\Ima \Psi(\xi))^2
     }
     \leq
     \frac{1}{\lambda+\Real \Psi(\xi)}
     .	
\end{align}
By the assumption, we can take constants $0<c_1<1< c_2<\infty$ and $R>0$ so that
\begin{align}
    c_1\leq  \frac{\Real \Psi(\xi)}{f(\xi)}\leq c_2
    \quad
    \text{and}
    \quad
    \frac{|\Ima \Psi(\xi)|}{f(\xi)}\leq c_2
    \quad
    \text{for any $\xi\geq R$.}	
\end{align}
Then for any $\lambda>0$ and $\xi\geq R$ we have
\begin{align}
    \Real \bigg(\frac{1}{\lambda+\Psi(\xi)}\bigg)\leq
    \frac{1}{\lambda+ c_1 f(\xi)}
    \leq
    \frac{1}{c_1}\frac{1}{\lambda+f(\xi)},	
\end{align}
and
\begin{align}
	\Real \bigg(\frac{1}{\lambda+\Psi(\xi)}\bigg)
	\geq
	\frac{\lambda+c_1 f(\xi)}{2(\lambda+c_2 f(\xi))^2}
	\geq
	\frac{c_1}{2c_2^2} \frac{1}{\lambda+f(\xi)}.
\end{align}
Therefore we can make the following estimate of $r_\lambda(0)$:
\begin{align}
	r_\lambda(0)
	&=\frac{1}{\pi}\int_{0}^\infty \Real \bigg(\frac{1}{\lambda+\Psi(\xi)}\bigg)
	\\
	&\leq
	\frac{1}{c_1\pi}\int_{R}^\infty\frac{d\xi}{\lambda+f(\xi)}+
	\frac{R}{\lambda}
	\sim
	\frac{1}{c_1\beta\sin(\pi/\beta)}
	\frac{f^{-1}(\lambda)}{\lambda}
	,
	\quad
	\text{as $\lambda\to\infty$}.
\end{align}
Here the last asymptotic relation can be shown in the same way as Lemma \ref{Lem:RV-levy-local}. 
Similarly,
\begin{align}
	r_\lambda(0)	&=\frac{1}{\pi}\int_{0}^\infty \Real \bigg(\frac{1}{\lambda+\Psi(\xi)}\bigg)
	\\
	&\geq
	\frac{c_1}{2 c_2^2 \pi}\int_R^\infty \frac{d\xi}{\lambda+f(\xi)}
	\sim
	\frac{c_1}{2 c_2^2\beta \sin(\pi/\beta)}\frac{f^{-1}(\lambda)}{\lambda},
	\quad
	\text{as $\lambda\to\infty$,} 
\end{align}
as desired.
\end{proof}

By the L\'evy--Khintchine formula, $\Psi$ can be represented as
\begin{align}
	\Psi(\xi)
	=
	i\xi \gamma +\frac{1}{2}\xi^2\sigma^2
	+
	\int_{\bR\setminus\{0\}}
	(1-e^{i\xi x}+i\xi x\mathbbm{1}_{\{|x|\leq  1\}})\,\nu(dx),
	\quad
	\xi\in\bR,
	\label{eq:LK:characteristic}
\end{align}
where $\gamma\in \bR$, $\sigma\geq0$ and a measure $\nu(dx)$ on $\bR\setminus\{0\}$ such that $\int_{\bR\setminus\{0\}}(1\wedge x^2)\,\nu(dx)<\infty$. The triplet $(\gamma, \sigma, \nu)$ is uniquely determined by $\Psi$. Conversely, for any triplet $(\gamma, \sigma, \nu)$, there exists a real-valued L\'evy process whose characteristic exponent $\Psi$ is given by \eqref{eq:LK:characteristic}.
The constant $\sigma$ and the measure $\nu$ are called the gaussian coefficient and the L\'evy measure, respectively. 
We can also give sufficient conditions for $\Psi$ in terms of the triplet $(\gamma,\sigma, \nu)$ to satisfy the assumptions of Lemmas \ref{Lem:RV-levy-local} and \ref{Lem:BRV-levy-local}.
Set
\begin{align}
    \bar{\nu}_{+}(x)=\nu(x,\infty)
    \quad
    \text{and}
    \quad
    \bar{\nu}_-(x)=\nu(-x,\infty),
\quad x>0.
\end{align}

\begin{Nota}
Let $f(x)$, $g(x)$ be functions defined on a neighborhood of $a\in[-\infty,\infty]$.
Assume $g(x)\neq 0$ on a further small neighborhood of $a$.
For $c\in\bC$, the notation $f(x)\sim c g(x)$ as $x\to a$ means $f(x)/g(x)\to c$ as $x\to a$ (even if $c=0$, following \cite[Preface, page xix]{BinGolTeu}).	
\end{Nota}

The following lemma is a corollary of the convergence theorem for infinitely divisible distributions \cite[Theorem 15.14]{Kal} and the L\'evy--Khintchine representation for stable distributions (see for example \cite[Section 1.2.6]{Kyp}). See also Ivanovs \cite[Theorem 2 (iii)]{Iva}. Nevertheless we give a more direct proof for the convenience of the readers. 

\begin{Lem}\label{Lem:RV-short-time-Levy-characteristic}
Let $\beta\in(1,2)$, $c_+, c_- \geq 0$ be constants with $c_+ + c_->0$ and let $\ell:(0,\infty)\to(0,\infty)$ be a function that is slowly varying at $0+$. Assume the following conditions are satisfied.
\begin{enumerate}[label=\upshape(\roman*)]
\item $\displaystyle \bar{\nu}_\pm(x)\sim c_{\pm} x^{-\beta}\ell(x)$, as $x\to0+$.
\item $\sigma=0$.

\end{enumerate}
Then  
\begin{align}
    \Psi(\xi)
    &\sim 
    \Gamma(1-\beta) \Big(c_+e^{-i\pi\beta/2}+c_-e^{i\pi\beta/2}\Big)\xi^\beta\ell(1/\xi),
    \\
    &=
    C
    \bigg(1-i\delta \tan \frac{\pi\beta}{2}\bigg)\xi^\beta\ell(1/\xi),
    \quad
    \text{as $\xi\to\infty$,}
    \label{RV-short-time-Levy-characteristic}
\end{align}
where 
\begin{align}
	C=\Gamma(1-\beta)(c_+ + c_-)\cos\frac{\pi\beta}{2}>0
	\quad\text{and}\quad
	\delta=\frac{c_+-c_-}{c_+ + c_-}\in[-1,1].
\end{align}
In particular, the condition \eqref{Asmp:Levy-integrability} is satisfied, and $|\Psi(\xi)|$ is regularly varying at $\infty$ with index $\beta$, and $\Psi(\xi)/|\Psi(\xi)|\to e^{i\theta}$ as $\xi\to \infty$ with
\begin{align}
\theta=\arctan\bigg(-\delta \tan \frac{\pi\beta}{2}\bigg)
    =
    \arctan \bigg(\delta \tan \bigg(1-\frac{\beta}{2}\bigg)\pi\bigg),
\end{align}
and hence $|\theta|\leq (1-\beta/2)\pi$. Here it is understood that $\arctan x$ takes values in $(-\pi/2, \pi/2)$. 
\end{Lem}

\begin{proof}
Let $\xi>1$. It is easily seen that
\begin{align}
    \Psi(\xi)
    &=\int_{0<|x|\leq 1}(1-e^{i\xi x}+i\xi x)\,\nu(dx) +O(\xi)
    \\
    &=\xi^2\int_0^\infty e^{i\xi x}\bigg(\int_{x\wedge 1}^1 \bar{\nu}_+(y)\,dy\bigg)\,dx
    \\
    &\phantom{=\,}+\xi^2\int_0^\infty e^{-i\xi x}
    \bigg(\int_{x\wedge 1}^1\bar{\nu}_-(y)\,dy\bigg)\,dx
    +O(\xi)
    ,
    \quad
    \text{as $\xi\to \infty$.}	
    \label{Lem:RV-short-time-Levy-characteristic-1}
\end{align}
The function $(0,\infty)\ni x \mapsto \int_{x\wedge 1}^1 \bar{\nu}_+(y)\,dy\in[0,\infty)$ is non-increasing, and
by the assumption,
\begin{align}
	\int_{x\wedge 1}^1 \bar{\nu}_+(y)\,dy
	\sim
	\frac{c_+}{\beta-1}x^{-\beta+1}\ell(x),
	\quad
	\text{as $x\to0+$}.
\end{align}
Hence we can use Abelian theorem for Fourier integrals (see for example \cite[Theorem 4.10.3]{BinGolTeu}) to get
\begin{align}
   \int_0^\infty e^{i\xi x}
   \bigg(\int_{x\wedge 1}^1 \bar{\nu}_+(y)\,dy\bigg)\,dx
    &=
    \int_0^\infty e^{ix}\bigg(\int_{(x/\xi)\wedge 1}^1 \bar{\nu}_+(y)\,dy\bigg)\,\frac{dx}{\xi}
    \\
    &\sim
    \bigg(\int_0^\infty e^{ix} x^{-1+(2-\beta)}\,dx\bigg) \frac{c_+}{\beta-1}\xi^{-2+\beta}\ell(1/\xi)
    \\
    &=
    \frac{\Gamma(2-\beta)}{\beta-1} (-i)^{-2+\beta} c_+ \xi^{-2+\beta}\ell(1/\xi)
    \\
    &=
   \Gamma(1-\beta)c_+ e^{-i\pi\beta/2}\xi^{-2+\beta} \ell(1/\xi),
   \quad
   \text{as $\xi\to\infty$.}
   \label{Lem:RV-short-time-Levy-characteristic-2}
\end{align}
Here we used
\begin{align}
    \int_0^\infty e^{-zx}x^{-1+\rho}\,dx
    =
    \Gamma(\rho)z^{-\rho}=\Gamma(\rho)e^{-\rho\log z},
    \label{eq:gamma-function}
\end{align}
for any $z\in\bC\setminus\{0\}$ with $\Real z\geq0$ and $0<\rho<1$. It is understood that the left-hand side of \eqref{eq:gamma-function} is the improper integral when $\Real z=0$. 
Similarly,
\begin{align}
	\int_0^\infty e^{-i\xi x}
    \bigg(\int_{x\wedge 1}^1\bar{\nu}_-(y)\,dy\bigg)\,dx
    \sim
	\Gamma(1-\beta)c_- e^{i\pi\beta/2}\xi^{-2+\beta} \ell(1/\xi),
	\quad
	\text{as $\xi\to\infty$.}
	\label{Lem:RV-short-time-Levy-characteristic-3}
\end{align}
Combining \eqref{Lem:RV-short-time-Levy-characteristic-1} with \eqref{Lem:RV-short-time-Levy-characteristic-2} and \eqref{Lem:RV-short-time-Levy-characteristic-3}, we obtain \eqref{RV-short-time-Levy-characteristic}. The rest of the proof is straightforward.
\end{proof}

We can also formulate a long-range version of the previous lemma.

\begin{Lem}\label{Lem:RV-long-time-Levy-characteristic}
Let $\beta\in(1,2)$, $c_+, c_- \geq 0$ be constants with $c_+ + c_->0$ and let $\ell:(0,\infty)\to(0,\infty)$ be a function that is slowly varying at $\infty$. Assume the following conditions are satisfied.
\begin{enumerate}[label=\upshape(\roman*)]

\item\label{Lem:RV-long-time-Levy-characteristic-1} $\displaystyle \bar{\nu}_\pm(x)\sim c_{\pm} x^{-\beta}\ell(x)$, as $x\to\infty$.
\item\label{Lem:RV-long-time-Levy-characteristic-2} 
$\gamma=\int_{|x|> 1} x\,\nu(dx)$.
\end{enumerate}
Then  
\begin{align}
    \Psi(\xi)
    &\sim 
    \Gamma(1-\beta) \Big(c_+e^{-i\pi\beta/2}+c_-e^{i\pi\beta/2}\Big)\xi^\beta\ell(1/\xi),
    \quad
    \text{as $\xi\to0+$.}
    \label{RV-long-time-Levy-characteristic}
\end{align}	
\end{Lem}

\begin{Rem}
The condition \ref{Lem:RV-long-time-Levy-characteristic-1} implies $\int_{|x|> 1}|x|\,\nu(dx)<\infty$. 
The conditions  \ref{Lem:RV-long-time-Levy-characteristic-1} and \ref{Lem:RV-long-time-Levy-characteristic-2} imply that $X_t$ is integrable and $\bE X_t=0$ for all $t\geq0$. See \cite[Example 25.12]{Sat}.	
\end{Rem}

\begin{proof}
Let $0<\xi<1$.
By the assumption,
\begin{align}
     \Psi(\xi)
     &=
     \int_{\bR\setminus\{0\}} (1-e^{i\xi x}+i\xi x)\, \nu(dx)
     +O(\xi^2)
      \\
    &=\xi^2\int_0^\infty e^{i\xi x}\bigg(\int_{x}^\infty \bar{\nu}_+(y)\,dy\bigg)\,dx
    \\
    &\phantom{=\,}+\xi^2\int_0^\infty e^{-i\xi x}
    \bigg(\int_{x}^\infty\bar{\nu}_-(y)\,dy\bigg)\,dx
    +O(\xi^2),
     \quad
     \text{as $\xi\to 0$+.}	
\end{align}
The rest of the proof is almost the same as that of Lemma \ref{Lem:RV-short-time-Levy-characteristic}. So we omit it.	
\end{proof}

Set
\begin{align}
 	\bar{\nu}(x)=\bar{\nu}_+(x)+\bar{\nu}_-(x)=\nu(\bR\setminus[-x,x]),
 	\quad
 	x>0.
\end{align}
The following lemma gives a sufficient condition for $\Psi$ to satisfy the assumptions in Lemma \ref{Lem:BRV-levy-local} in the regime $\xi\to\infty$.
 
\begin{Lem}\label{Lem:BRV-short-time-Levy-characteristic}
Let $\beta\in(1,2)$, $c_+, c_- \geq 0$ be constants with $c_+ + c_->0$ and let $\ell:(0,\infty)\to(0,\infty)$ be a function that is slowly varying at $0+$. Assume the following conditions are satisfied.
\begin{enumerate}[label=\upshape(\roman*)]

\item $\bar{\nu}(x)\asymp x^{-\beta}\ell(x)$, as $x\to0+$.	
\item $\sigma=0$.
\end{enumerate}
Then
\begin{align}
    \Real \Psi(\xi)\asymp \xi^\beta\ell(1/\xi)
    \quad\text{and}\quad
    \Ima \Psi(\xi)=O\Big(|\xi|^\beta \ell(1/\xi)\Big),
    \quad
    \text{as $\xi\to\infty$,}	
\end{align}
and hence the condition \eqref{Asmp:Levy-integrability} is also satisfied.
\end{Lem}
\begin{proof} By the assumption, there exist constants $0<C_1\leq C_2<\infty$ and $R>1$ such that
\begin{align}
	C_1\leq \frac{\bar{\nu}(x)}{x^{-\beta}\ell(x)},
	\,\,
	\frac{\int_x^1\bar{\nu}(y)\,dy}{x^{-\beta+1}\ell(x)},
	\,\,
	\frac{\int_0^x y^n \bar{\nu}(y)\,dy}{x^{-\beta+n+1}\ell(x)}
	\leq
	C_2,
	\quad
	0<x<\frac{1}{R},
	\,\, n=1,2.
\end{align}
Note that
\begin{align}
    \frac{1}{3}x^2\mathbbm{1}_{\{|x|\leq 1\}}\leq \bigg(\frac{1}{2}x^2-\frac{1}{24}x^4\bigg)\mathbbm{1}_{\{|x|\leq 1\}}\leq 1-\cos x\leq \bigg(\frac{1}{2}x^2\bigg)\wedge 2,
    \quad
    x\in\bR.
\end{align}
Hence for $\xi>R$,
\begin{align}
	 \Real \Psi(\xi)
	 &=\int_{\bR\setminus\{0\}} (1-\cos (\xi x))\,\nu(dx)
	 \\
	 &\geq
	 \frac{\xi^2}{3}\int_{0<|x|\leq 1/\xi} x^2\, \nu(dx)
	 =
	 \frac{\xi^2}{3}\int_0^{1/\xi} y\,\bar{\nu}(y)\,dy
	 \geq
	 \frac{C_1}{3} \xi^\beta \ell(1/\xi),
\end{align}
and
\begin{align}
	\Real \Psi(\xi)
	&\leq
	\frac{\xi^2}{2}\int_{|x|\leq 1/\xi}
	x^2\,\nu(dx)+
	2\bar{\nu}(1/\xi)
	\leq
	\frac{5C_2}{2} \xi^\beta \ell(1/\xi).
\end{align}
Therefore $\Real \Psi(\xi)\asymp \xi^\beta \ell(1/\xi)$, as $\xi\to\infty$. Similarly, we have
\begin{align}
    |x-\sin x| \leq \bigg(\frac{1}{6}|x|^3\bigg)\wedge (2|x|),
    \quad
    x\in\bR,
\end{align}
and hence for $\xi>R$,
\begin{align}
    |\Ima \Psi(\xi)|
    &\leq
    \int_{0<|x|\leq 1}
    |\xi x-\sin(\xi x)|\,\nu(dx) +O(\xi)
    \\
    &\leq
    \frac{\xi^3}{6}\int_{0<|x|\leq 1/\xi}
    |x|^3\,\nu(dx)
    +2\xi \int_{1/\xi<|x|\leq 1}|x|\,\nu(dx)+O(\xi)
    \\
    &=
    \frac{\xi^3}{2}\int_0^{1/\xi} x^2\,\bar{\nu}(x)\,dx
    +
    2\xi\int_{0}^1 
    \bigg(\bar{\nu} \bigg(\frac{1}{\xi}\vee x\bigg)- \bar{\nu}(1)\bigg)\,dx
    +
    O(\xi)
    \\
    &=
    \frac{\xi^3}{2}\int_0^{1/\xi} x^2\,\bar{\nu}(x)\,dx
    +
    2\bar{\nu}(1/\xi)+2\xi\int_{1/\xi}^1 \bar{\nu}(x)\,dx
    +O(\xi)
    \\
    &= O( \xi^\beta \ell(1/\xi)),
    \quad
    \text{as $\xi\to\infty$,}
    \end{align}
as desired.
\end{proof}

We can also consider the case $\xi\to0+$:

\begin{Lem}
Let $\beta\in(1,2)$, $c_+, c_- \geq 0$ be constants with $c_+ + c_->0$ and let $\ell:(0,\infty)\to(0,\infty)$ be a function that is slowly varying at $\infty$. Assume the following conditions are satisfied.
\begin{enumerate}[label=\upshape(\roman*)]

\item $\bar{\nu}(x)\asymp x^{-\beta}\ell(x)$, as $x\to\infty$.	
 \item $\gamma=\int_{|x|> 1} x\,\nu(dx)$.
\end{enumerate}
Then
\begin{align}
    \Real \Psi(\xi)\asymp \xi^\beta\ell(1/\xi)
    \quad\text{and}\quad
    \Ima \Psi(\xi)=O\Big(|\xi|^\beta \ell(1/\xi)\Big),
    \quad
    \text{as $\xi\to0+$.}
\end{align}	
\end{Lem}

\begin{proof}
Let $\xi>0$.
Under the assumptions, we have
\begin{align}
    \Real \Psi(\xi)
    =
    \int_{\bR\setminus\{0\}}
    (1-\cos\xi x)\,\nu(dx)+O(\xi^2),
    \quad
    \text{as $\xi\to0+$,}
\end{align}
and
\begin{align}
     |\Ima \Psi(\xi)|
     &=
     \bigg|\int_{\bR\setminus\{0\}}(\xi x-\sin \xi x)\,\nu(dx)\bigg|
     \\
     &\leq
     \frac{\xi^3}{6}\int_{0<|x|\leq 1/\xi} |x|^3\,\nu(dx)
     +
     2\xi\int_{1/\xi <|x|}|x|\,\nu(dx)
     \\
     &=
     \frac{\xi^3}{2}\int_0^{1/\xi} y^2\,\bar{\nu}(y)\,dy
     +
     2\bar{\nu}(1/\xi)+2\xi\int_{1/\xi}^\infty \bar{\nu}(y)\,dy.
\end{align}
The rest of the proof is almost the same as that of Lemma \ref{Lem:BRV-short-time-Levy-characteristic}, so we omit it.
\end{proof}

\section{Examples: Brownian motion and Stable Processes on Unbounded Sierpi\'nski Gasket}\label{Sec:SG}

Let $X$ be a Hunt process satisfying the conditions \ref{Asmp-Hunt}, \ref{Asmp-Dual}, \ref{Asmp-u1}, \ref{Asmp-Dual-P}, \ref{Asmp-sym}, \ref{Asmp-sym-strong} of Section \ref{Sec:Dual}. By the results in Section \ref{Sec:Dual}, the inverse local time $(\tau, \bP_a)$ is a complete subordinator with Laplace exponent $\Phi(\lambda)=1/r_\lambda(a,a)$ and potential measure $U(dt)=p_t(a,a)\,dt$.
Therefore we can apply Theorems \ref{Thm:LV-LDE-long} (2) and \ref{Thm:LV-LDE-short} (2) if $r_\lambda(a,a)$ is bounded above and below by a regularly varying function with index $-\alpha$ for some $\alpha\in(0,1)$. Such on-diagonal estimates hold when $X$ satisfies suitable sub-Gaussian or stable-like estimates. In the following we focus on these estimates for the Brownian motion and stable processes on the unbounded Sierpi\'nski gasket, and then apply our results to these processes.

\begin{Ex}[Brownian motion on the unbounded Sierpi\'nski gasket]\label{Ex:BM-on-SG-revisited}
We use the same notation as in Example \ref{Ex:SG}.
Let $X$ be the Brownian motion on the unbounded Sierpi\'nski gasket $\hat{K}$.
Set 
\begin{align}
	d_f=\frac{\log3}{\log 2},
	\quad
	d_s=\frac{2\log3}{\log 5}
	\quad
	\text{and}
	\quad
	d_w=
    \frac{2d_f}{d_s}
    =
    \frac{\log 5}{\log 2}.
\end{align}
The constants $d_f$, $d_s$ and $d_w$ are called the fractal (or Hausdorff) dimension, the spectral dimension and the walk dimension of $\hat{K}$, respectively.
By \cite[Theorems 1.4 and 1.5]{BarPer},  
$X$ is a Feller process, and hence a Hunt process on $(\hat{K}, m)$ with $L^2(m)$-strongly continuous semigroup. Moreover, it admits a  transition density $p_t(x,y)=p_t(y,x)$ that is jointly continuous with respect to $(t,x,y)$.
Furthermore, the following sub-Gaussian heat kernel estimate holds: there exist constants $C_1,C_2,C_3,C_4\in(0,\infty)$ such that
\begin{align}
	\frac{C_1}{t^{d_s/2}}
	&\exp\bigg(\!-C_2\bigg(\frac{|x-y|^{d_w}} {t}\bigg)^{1/(d_w-1)}\bigg)
	\\
	&\leq
	p_t(x,y)
	\leq 
	\frac{C_3}{t^{d_s/2}}\exp\bigg(\!-C_4\bigg(\frac{|x-y|^{d_w}} {t}\bigg)^{1/(d_w-1)}\bigg),
	\label{Ex:SG-sub-Gaussian}
\end{align}
for all $(t,x,y)\in (0,\infty)\times \hat{K}\times \hat{K}$.
In particular, we have \eqref{Ex:SG-HK-diagonal} and \eqref{Ex:SG-PD-diagonal}. The conditions   \ref{Asmp-Hunt}, \ref{Asmp-Dual}, \ref{Asmp-u1}, \ref{Asmp-Dual-P}, \ref{Asmp-sym}, \ref{Asmp-sym-strong} of Section \ref{Sec:Dual} are satisfied for any $a\in\hat{K}$. Therefore we can apply Theorems \ref{Thm:LV-LDE-long} (2) and \ref{Thm:LV-LDE-short} (2) with $\alpha=1-(d_s/2)\in(0,1)$ to this setting. 
 
\end{Ex}

\begin{Rem}[Gaussian and sub-Gaussian heat kernel estimates]
Gaussian and sub-Gaussian heat kernel estimates such as \eqref{Ex:SG-sub-Gaussian} have been extensively studied in  geometric analysis and probability theory. For example, Li and Yau \cite{LiYau} established Gaussian heat kernel estimates on complete Riemannian manifolds with nonnegative Ricci curvature. Subsequent work by Grigor'yan \cite{Gri91} and Saloff-Coste \cite{Sal92} developed Gaussian heat kernel estimates under assumptions of various analytic and geometric inequalities. 
Barlow and Perkins \cite{BarPer} demonstrated the subdiffusive behavior of Brownian motion on the Sierpi\'nski gasket by proving the sub-Gaussian estimate \eqref{Ex:SG-sub-Gaussian}. Their result was later extended to various classes of fractal diffusions, including Brownian motions on the Sierpi\'nski carpet (Barlow and Bass \cite{BarBas92}), nested fractals (Kumagai \cite{Kum93}) and  affine nested fractals (Fitzsimmons, Hambly,  and Kumagai \cite{FitHamKum}). See also Barlow \cite{Bar} and Kigami \cite{Kig01,Kig12}. We further refer the reader to Kajino \cite{Kaj13b,Kaj13a} for results on on-diagonal oscillations of heat kernels for such processes, and to 
Grigor'yan and Telcs \cite{GriTel12}, 
Barlow, Grigor'yan, and Kumagai \cite{BarGriKum}, Kajino and Murugan \cite[Section 4]{KajMur23}, Eriksson-Bique \cite{EriBiq} and the references therein for relations between sub-Gaussian heat kernel estimates and various analytic and geometric inequalities.
\end{Rem}

\begin{Ex}[Subordinate Brownian motion on the unbounded Sierpi\'nski gasket]
 Under the setting of Example \ref{Ex:BM-on-SG-revisited}, let $\beta\in(0,1)$ and let $Y=(Y(t))_{t\geq0}$ be a $\beta$-stable subordinator independent of $X$ with Laplace transform $\bE\exp(-\lambda Y(t))=\exp(-\lambda^{\beta} t)$, $\lambda,t\geq0$. 
Let us consider the subordinate process $Z=(X_{Y(t)})_{t\geq0}$. Since $X$ is a Feller process, $Z$ is also a Feller (and hence Hunt) process. See for example \cite[Theorem 32.1]{Sat}.
 As shown in Bogdan, St\'os, and Sztonyk \cite[Theorem 3.1]{BogStoSzt-b}, $Z$ admits a jointly continuous transition probability density $p^{(\beta)}_t(x,y)=p^{(\beta)}_t(y,x)
 =\int_0^\infty p_s(x,y)\,\bP(Y(t)\in ds)$
  satisfying  the stable-like estimate
 \begin{align}
 	p_t^{(\beta)}(x,y)
 	\asymp
 	\frac{1}{t^{d_f/(\beta d_w)}}
 	\wedge
 	\frac{t}{|x-y|^{d_f+\beta d_w}},
 	\quad
 	\text{$(t,x,y)\in (0,\infty)\times \hat{K}\times \hat{K}$}.
 	\label{Ex:SG-stable-like}
 \end{align}
In particular,  
\begin{align}
	p_t^{(\beta)}(x,x)\asymp t^{-d_f/(\beta d_w)}= t^{-d_s/(2\beta)},
	\quad
	\text{
 	$(t,x)\in (0,\infty)\times \hat{K}$}. 
\end{align}
Therefore, all the conditions in Section \ref{Sec:Dual} other than \ref{Asmp-u1} are satisfied.
Denote by 
\begin{align}
r_\lambda^{(\beta)}(x,y)=\int_0^\infty e^{-\lambda t} p_t^{(\beta)}(x,y)\,dt,
\quad (\lambda, x,y)\in (0,\infty)\times \hat{K}\times \hat{K},
\end{align}
  the $\lambda$-resolvent density of $Z$. 
When $0<\beta\leq d_s/2$, we have $r^{(\beta)}_\lambda(x,x)=\infty$, which violates the condition \ref{Asmp-u1}, and hence we cannot apply our results.
In the following we focus on the case $d_s/2<\beta<1$. 
Then $r_\lambda^{(\beta)}(x,y)$ satisfies the on-diagonal estimate
\begin{align}
	r^{(\beta)}_\lambda(x,x)
	\asymp
	\lambda^{-1+(d_s/(2\beta))},
	\quad
	\text{$(\lambda, x)\in(0,\infty)\times \hat{K}$}.
\end{align}
Hence the process $Z$ satisfies the condition \ref{Asmp-u1} for any $a\in\hat{K}$. Therefore we can apply Theorems \ref{Thm:LV-LDE-long} (2) and \ref{Thm:LV-LDE-short} (2) with $\alpha=1-(d_s/(2\beta))\in(0,1-(d_s/2))$ and $d_s/2<\beta<1$ to this setting.
\end{Ex}

\begin{Rem}[Stable-like heat kernel estimates]
Stable-like heat kernel estimates such as \eqref{Ex:SG-stable-like} are analogous to the classical estimates for transition densities of symmetric stable processes on the Euclidean space $\bR^n$, whose study goes back to Blumenthal and Getoor \cite{BluGet60}. 
Bass and Levin \cite{BasLev} established stable-like estimates of transition probabilities for symmetric jump Markov chains on the integer lattice $\bZ^d$.
Bogdan, St\'os, and Sztonyk \cite{BogStoSzt-b} demonstrated heat kernel estimates for stable processes on $d$-sets. These processes are  obtained by applying
stable subordination to diffusions whose heat kernels satisfy sub-Gaussian
estimates.
In the setting of Example \ref{Ex:SG-stable-like}, the process $Z$ is associated with the symmetric Dirichlet form $(\mathcal{E}^{(\beta)},\mathcal{F}^{(\beta)})$ and  the jump kernel $J^{(\beta)}(x,y)$ in the form
\begin{align}
    J^{(\beta)}(x,y)&=
    \frac{1}{2}\int_0^\infty p_s(x,y)\frac{\beta s^{-1-\beta}}{\Gamma(1-\beta)}\,ds,
    \quad \text{$(x,y)\in \hat{K}\times \hat{K}$ with $x\neq y$,}
    \\
	\mathcal{F}^{(\beta)}&=\bigg\{u\in L^2(m)\::\: \int_{\hat{K}\times \hat{K}}
	(u(x)-u(y))^2 J^{(\beta)}(x,y)\,m(dx)m(dy)<\infty
	\bigg\},
	\\
	\mathcal{E}^{(\beta)}(u,v)
	&=
	\int_{\hat{K}\times \hat{K}}
	(u(x)-u(y))(v(x)-v(y)) J^{(\beta)}(x,y)\,m(dx)m(dy),
	\quad
	u,v\in \cF^{(\beta)}.
\end{align}
Here $\beta s^{-1-\beta}/\Gamma(1-\beta)$ denotes the density of the L\'evy measure of  the $\beta$-stable subordinator $Y$. For the general formula for the Dirichlet form of a subordinate process, see \^{O}kura \cite[Theorem 2.1]{Oku}.  It is  easily seen that
\begin{align}
      J^{(\beta)}(x,y)	
      \asymp
      |x-y|^{-(d_f+\beta d_w)},
      \quad
      \text{for all $(x,y)\in \hat{K}\times \hat{K}$ with $x\neq y$.}
\end{align}
See for example the proof of \cite[Corollary 6.2]{BogStoSzt-b}.
Chen and Kumagai \cite{CheKum03} established stable-like heat kernel estimates for symmetric pure jump processes on $d$-sets associated with non-local Dirichlet forms whose jump kernels are bounded above and below by $|x-y|^{-(d+\gamma)}$ with $\gamma\in(0,2)$.
They subsequently extended these results to mixed-type jump kernels in \cite{CheKum08}.
Taking $d=d_f$ and $\gamma=\beta d_w$,
their results imply \eqref{Ex:SG-stable-like} whenever $\gamma\in(0,2)$, or equivalently, $0<\beta <2/d_w$. The restriction of the parameter $\gamma$ was later removed by Chen, Kumagai, and Wang \cite{CheKumWan21} under additional assumptions. More precisely, they established equivalent conditions for stable-like heat kernel estimates in a general setting.
See also Murugan \cite{Mur26+} for a recent progress on  characterizations of stable-like heat kernel estimates.
\end{Rem}

%
%
%
%
%
%
%
%
%

\subsection*{Acknowledgements}

TM was supported by JSPS Grants-in-Aid for Scientific Research (KAKENHI) Grant Numbers 
JP23K25773 and  JP25K17266.
KN was supported by JSPS Grants-in-Aid for Scientific Research (KAKENHI) Grant Number JP21K13807.
TS was supported by JSPS  Grants-in-Aid for Scientific Research (KAKENHI) Grant Number JP24K16948.
This work was also supported by the ISM Cooperative Research Symposium (2025-ISMCRP-5007).
During the preparation of this manuscript, the authors used ChatGPT and Google Gemini to obtain suggestions for improving the presentation and exposition of the results. All suggestions generated by these tools were critically reviewed and independently verified by the authors, who take full responsibility for the final content of the manuscript.

\bibliographystyle{plain}

\end{document}